\documentclass[10pt]{article}

\usepackage{a4wide}
\usepackage{amssymb}
\usepackage{amsfonts}
\usepackage{amsmath}
\input xy
\xyoption{arrow} \xyoption{matrix}

\date{}

\newtheorem{proposition}{Proposition}[section]
\newtheorem{theorem}[proposition]{Theorem}
\newtheorem{lemma}[proposition]{Lemma}

\newtheorem{corollary}[proposition]{Corollary}

\def\der{\partial }

\def\nFM0{{\nu }_{F,M_0}}
\def\nFN0{{\nu }_{F,N_0}}
\def\nGN0{{\nu }_{G,N_0}}

\def\N0{ {\bf N}_0 }

\def\g{\gamma}
\def\v{\varphi}
\def\ra{\rightarrow}

\def\Xpm{X^{\pm }}

\def\s{\sigma}
\def\Z{{\bf Z }}

\def\l1{{\lambda}_1}

\def\a{\alpha}
\def\a0{ {\alpha }_0}
\def\a1{ {\alpha }_1}

\def\l{\lambda}

\def\nFGM0{{\nu }_{F,G,M_0}}

\def\nFN0{{\nu}_{F,N_0}}

\def\sm{{\sigma}^m}

\def\sm1{{\sigma}^{-1}}

\def\smtp1{{\sigma}^{-t+1}}

\def\S1{S^{-1}}

\def\Xpm1{X^{\pm 1}_1}

\def\sPM1{{\sigma }^{\pm 1}}
\def\sMP1{{\sigma }^{\mp 1 }}

\def\d{\delta}

\def\di{{\rm d.ind}}

\def\L{\Lambda}

\def\G{\Gamma}

\def\OO{{\cal O}}
\def\CA{{\cal A}}

\def\CD{{\cal D}}

\def\Ytm1{Y^{t-1}}
\def\Yim1{Y^{i-1}}

\def\Aut{{\rm Aut}}

\def\bA{\overline{A}}

\def\dim{{\rm dim }}

\def\Iso{{\rm Iso }}

\def\ker{ {\rm ker } }

\def\gcd{ {\rm gcd } }
\def\D{ \Delta }

\def\SL2Z{ {\rm SL}_2({\bf Z}) }

\def\th{ \theta }

\def\Gp1{ G^{1 , 1 } }
\def\P11{ P^{-1 , 1 } }
\def\Pp1{ P^{1 , 1 } }

\def\th{\theta}

\def\nCLsr{{}^\nu\kern-2pt {\cal L}^{\sigma , \rho  }}
\def\nP{{}^\nu \kern-2pt P}
\def\nL{{}^\nu\kern-2pt L}
\def\nLL{{}^\nu\kern-2pt \Lambda}
\def\nPsr{{}^\nu\kern-2pt P^{\sigma , \rho  }}
\def\nLsr{{}^\nu\kern-2pt L^{\sigma , \rho  }}
\def\nuCL{{}^\nu\kern-2pt  {\cal L}}
\def\nCLsr{{}^\nu\kern-2pt {\cal L}^{\sigma , \rho  }}
\def\nCL1m{{}^\nu\kern-2pt {\cal L}^{-1 , 1  }}
\def\x1nu{x^\frac{1}{\nu}}
\def\xm1nu{x^{-\frac{1}{\nu}}}

\def\rad{{\rm rad}}

\def\ra{\rightarrow }

\def\CB{{\cal B}}

\def\nAM0{{\nu }_{{\cal A},M_0}}
\def\nAN0{{\nu }_{{\cal A},N_0}}

\def\hCA{ \widehat{\CA} }

\def\bp{\overline{p}}

\def\ga{\mathfrak{a}}

\def\gc{\mathfrak{c}}

\def\gn{\mathfrak{n}}

\def\di!{\frac{\der^i}{i!}}
\def\dik!{\frac{\der^k_i}{k!}}

\def\N{\mathbb{N}}

\def\0{\overline{0}}
\def\1{\overline{1}}

\def\Ln1{\L_{n,\overline{1}}}

\def\a1{a_{\overline{1}}}

\def\S{\Sigma}

\def\vn1{\overrightarrow{n-1}}

\def\Sh{{\rm Sh}}

\def\Q{\mathbb{Q}}

\def\mA{\mathbb{A}}

\def\mL{\mathbb{L}}

\def\Sub{{\rm Sub}}

\def\mS{\mathbb{S}}
\def\mJ{\mathbb{J}}
\def\mI{\mathbb{I}}

\def\ann{{\rm ann}}

\def\mM{\mathbb{M}}
\def\mT{\mathbb{T}}
\def\ind{{\rm ind}}

\def\mU{\mathbb{U}}

\def\K1{{\rm K}_1}

\def\hmI1{\widehat{\mI_1}}
\def\tmI1{\widetilde{\mI_1}}
\def\tmJ1{\widetilde{\mJ_1}}
\def\hB1{\widehat{B_1}}
\def\hCB1{\widehat{\CB_1}}

\def\ga{\mathfrak{a}}

\def\Z{\mathbb{Z}}

\def\mL{\mathbb{L}}

\def\mB{\mathbb{B}}

\def\CB{{\cal B}}

\def\hB{\hat{B}}

\def\mC{\mathbb{C}}

\def\dec{{\rm dec}}
\def\Rel{{\rm Rel}}
\def\mO{\mathbb{O}}

\def\Amon{A_{\rm mon}}
\def\bAmon{\overline{A}_{\rm mon}}

\begin{document}

\author{V. V.   Bavula 
} 
 
\title{The Question of Arnold on classification of co-artin subalgebras in singularity theory}

\maketitle 
\begin{abstract}

In \cite[Section 5, p.32]{Arnold-1998}, Arnold writes: ``Classification of singularities of curves can be interpreted in dual terms as a description of `co-artin' subalgebras of finite co-dimension in the algebra of formal series in a single variable (up to isomorphism of the   algebra of formal series).'' In the paper,  such a description is obtained  but up to isomorphism of algebraic curves (i.e. this description is finer).

Let $K$ be an algebraically closed field of arbitrary characteristic. The aim of the paper is to give a classification (up to isomorphism) of the set of subalgebras $\CA$ of the polynomial algebra $K[x]$ that contains the ideal $x^mK[x]$ for some $m\geq 1$. It is proven that the  set $\CA = \coprod_{m, \G }\CA (m, \G )$ is a disjoint union of affine algebraic varieties (where $\G \coprod \{0, m, m+1, \ldots \}$ is the semigroup of the singularity and $m-1$ is the Frobenius number). It is proven that each set $\CA (m, \G )$ is an affine algebraic variety and explicit generators and defining relations are given for the  algebra  of regular functions on $\CA (m ,\G )$.  An isomorphism criterion is given for the algebras in $\CA$. For each algebra $A\in \CA (m, \G)$,  explicit sets of generators and defining relations are given and the automorphism group $\Aut_K(A)$ is explicitly described.  The automorphism group of the algebra $A$ is finite iff the algebra $A$ is not isomorphic to a monomial algebra, and in this case $|\Aut_K(A)|<\dim_K(A/\gc_A)$ where $\gc_A$ is the conductor of $A$. The set of orders of the automorphism groups of the algebras in $\CA (m , \G )$ is explicitly described. \\


 {\em Mathematics subject classification
2010:  14H20, 14H37, 14R05, 14H10, 14J10. 

$${\bf Contents}$$
\begin{enumerate}
\item Introduction.
\item The canonical basis, generators and defining relations of algebras $\bA$ and $A$ where $A\in \CA (m, \G )$.
\item  Isomorphism problems and an explicit description of the automorphism groups of algebras   $A\in \CA (m)$.
\item Generators and defining relations for the algebra  $\OO (\CA (m, \G ))$ of regular functions on the algebraic variety $\CA (m , \G )$.

\end{enumerate}
}

\end{abstract}


\section{Introduction}

{\em Motivation:} A {\em singularity of a curve} means the germ of a holomorphic map of the complex line into complex space at a singular point, \cite{Arnold-1998}. Classifying curve singularities up to diffeomorphism is a classical problem in the theory of algebraic curves. 
 A singularity is called {\em simple} if all the singularities of the neighbouring mappings belong
to a finite set of equivalence classes.
Simple singularities of plane curves were classfied by J. W. Bruce and T. J. Gaffney \cite{Bruce-Gafney-82}, and simple singularities of space curves by C. G. Gibson and C. A. Hobbs \cite{Gibson-Hobbs-93}. The classification of simple singularities of curves is
described by Arnold  \cite{Arnold-1998}.
 
In \cite[Section 5, p.32]{Arnold-1998}, Arnold writes: ``Classification of singularities of curves can be interpreted in dual terms as a description of `co-artin' subalgebras of finite co-dimension in the algebra of formal series in a single variable (up to isomorphism of the   algebra of formal series),'' see the end of the Introduction for more details.   Such a description is obtained in the present paper but up to isomorphism of algebraic curves (i.e. this description is finer). 

The following notation will remain fixed throughout the paper: $K$ is an arbitrary algebraically closed field of arbitrary characteristic (many results of the paper are true for an arbitrary field); $K^\times := K\backslash \{ 0\}$;  algebra means $K$-algebra; $K[x]$ is a polynomial algebra in the variable  $x$ over $K$; $m\geq 2$ is a natural number; $\CA (m)$ is the set of all $K$-subalgebras $A$ of $K[x]$ such that $x^mK[x]\subset A$ but $x^{m-1}\not\in A$, i.e. the ideal 
 $x^mK[x]$ of $A$ is the largest ideal of $A$ which is also an ideal of $K[x]$, i.e. the ideal $x^mK[x]=\ann_A(K[x]/A)$ is the {\em conductor} of $A$ as the polynomial algebra $K[x]$ is the normal closure of $A$;
 $$\bA :=A/x^mK[x]\subseteq F:=F_m:=K[x]/x^mK[x];$$  
$\Aut_K(A)$ is the automorphism group of the $K$-algebra $A$. \\





{\bf The canonical basis of the algebra $\bA$.} Let $\mS (m)$ be the set that contains the empty set and all non-empty subsets $\G$ of the set $\{ 2,\ldots , m-2\}$ such that $\G +\G \subseteq \G \cup [m, \infty )$. Notice that, by definition, $m-1\not\in \G$; $\mS (m)\neq \emptyset$ iff $m\geq 4$. For $m=1,2,3$, the set $\CA (m)$ contains the only algebra $K+x^mK[x]$. So, we will assume that $m\geq 4$. 

Let $A\in \CA (m)$. The the ideal $(x)$ of the finite dimensional algebra $F=F_m$ is its radical. The $(x)$-adic filtration on $F$ is also the radical filtration. It induces a  filtration on the  subalgebra $\bA$ of $F=F_m$ and the associated graded algebra
(see (\ref{grbA}) and (\ref{grbA1}))
$$
{\rm gr}(\bA )=K\oplus \bigoplus_{\g \in \G_A}Kx^\g ,\;\; {\rm 
where}\;\; \G_A:=\{ \g \, | \, x^\g \in {\rm gr}(\bA ), 1\leq \g \leq m-1\}\in \mS (m),$$
is a graded algebra. For each $\G \in \mS (m)$, let $\CA (m, \G ) :=\{ A\in \CA (m)\, | \, \G_A=\G\}$. Then 
$$\CA (m)=\coprod_{\G\in \mS (m)}\CA (m,\G).$$
We will see that each set $\CA (m, \G )$ is an affine algebraic variety (Lemma \ref{a21Jun20}, Theorem \ref{XX12May20} and Theorem \ref{XY12May20}). In the literature, for each algebra $A\in \CA (m , \G )$, the subsemigroup  $\G \cup \{ 0, m, m+1, \ldots \}$ of the semigroup of natural numbers $(\N , +)$ is called the 
{\bf semigroup of the singularity} and the number $m-1$ is called the {\bf Frobenius number}.

Proposition \ref{A4May20} states that for 
 each algebra $A\in \CA (m, \G \}$, there is a {\em unique basis} $\{ 1, f_\g, \, | \, \g \in \G \}$ of the algebra $\bA$ such that  $$f_\g =x^\g +\sum_{\d \in C\G (\g )}\l_{\g \d}x^\d\;\; {\rm  where}\;\; \l_{\g \d}\in K\;\; {\rm and}\;\;  C\G (\g ) :=\{ \d \, | \, \d\not\in \G, \g <\d \leq m-1\}.$$
The basis $\{ 1, f_\g, \, | \, \g \in \G \}$ is called the {\bf canonical basis} of the algebra $\bA$. The canonical  basis has  many remarkable properties. It is used in many proofs of  this paper. In particular, the elements of the canonical basis are common eigenvectors for the automorphism group of the algebra $A$ (each element $\s\in \Aut_K(A)$ preserves the conductor of $A$ and hence induces an automorphism of the algebra $\bA$) and this is the key fact in finding the automorphism group $\Aut_K(A)$.\\

{\bf Generators and defining relations of the algebras $A\in \CA (m , \G )$ and $\bA$.}\\

Each non-empty set $\G \in \mS (m)$ is a  disjoint union $\G = \ind (\G ) \coprod \dec (\G )$ where $\ind (\G )$ is the set of {\em indecomposable} elements of $\G$ and $\dec (\G )$ is the set of {\em decomposable} elements of $\G$. An element $\g \in \G$ is called {\em decomposable} if $\g =\g_1+\g_2$ for some elements $\g_1, \g_2\in \G$. Every decomposable element is a linear combination of indecomposable ones with natural coefficients (which is highly non-unique, in general). We denote by $\dec (\G )_{\geq 2}$   
 the set of all decomposable elements that admits  at least 2 distinct linear presentations. Let $\ind (\G ) = \{ \nu_1, \ldots , \nu_s\}$ where $s=|\ind (\G )|$.  For each element $\g \in \G$, we fix a vector $a(\g )=( a(\g )_1, \ldots , a(\g )_s)\in \N^s$ such that $\g =\sum_{i=1}^sa(\g )_i\nu_i$.\\
  
The algebras in $\CA (m, \G)$ are partitioned into three isomorphism invariant classes: 

(I)  $|\ind (\G )|=1$, 

(II)  $|\ind (\G )|\geq 2$ and $\dec (\G )_{\geq 2}=\emptyset$,

(III) $|\ind (\G )|\geq 2$ and $\dec (\G )_{\geq 2}\neq \emptyset$.\\

For each algebra in the classes (I), (II) and (III) explicit sets of  generators and defining relations  are given in Lemma \ref{Aa7May20},  Corollary  \ref{Ab7May20} and  Theorem \ref{XA6May20}, respectively.\\

{\bf Isomorphism problems and the automorphism groups of algebras $A$.}\\

Theorem \ref{3May20}.(2) and Proposition \ref{a20Jul20} are isomorphism criteria for algebras in $\CA :=\coprod_{m\geq 2}\CA (m)$.  For algebras  $A, A'\in \CA$ such that $A\simeq A'$, Theorem \ref{3May20}.(3) describes the set $\Iso_K(A,A')$ of all $K$-algebra isomorphisms from  $A$ to $A'$.  Theorem \ref{28Apr20}.(1) is an explicit description of the automorphism groups of algebras in $\CA (m)$. Theorem \ref{28Apr20}.(2-4)  gives criteria for the automorphism group  to be infinite/finite/trivial (see also Corollary \ref{a2May20}).

Let $\mO (m)$ be the set of orders of all {\em finite} automorphism groups of algebras in $\CA (m)$. By Corollary \ref{a30Apr20}.(1),
$$ \mO (m) \subseteq \begin{cases}
\{ 1, \ldots , m-3\}& \text{if } 4\leq m \text{ is even},\\
\{ 1, \ldots , m-4\}& \text{if } 5\leq m \text{ is odd}.\\
\end{cases}
$$
Let $p$ be a prime number. A natural number $i\geq 1$ is a unique product $i=p^di_p$ for some natural numbers $d$ and $i_p$ such that $(p,i_p)=1$. The natural number $i_p$ is called the $p$-co-prime divisor of $i$. For a natural number $m\geq 4$, let   $u(p,m):=\max\{i_p\, | \, 1\leq i\leq m-3$, $i+j\leq m-1$ for some natural number $j=j(i)\geq 2$ such that $j\nmid m-1\}$.  By Corollary \ref{a30Apr20}.(1),
$$ \max \mO (m) =\begin{cases}
 m-3& \text{if } 4\leq m \text{ is even, char\,}K=0,\\
 m-4& \text{if } 5\leq m \text{ is odd, char\,}K=0,\\
  u(p,m)& \text{if } 4\leq m, \text{ char\,}K=p>0.\\
\end{cases}
$$
Theorem \ref{2May20} and Corollary \ref{c3May20}
 are explicit descriptions of the set $\mO (m)$. The set $\mO (m)$ is an intricate set (see (\ref{Lm1}),  (\ref{Lm2}) and (\ref{Lm3})). Lemma \ref{a23Apr20} and Lemma \ref{d3May20} are examples of algebras in $\CA (m)$ with sophisticated orders of automorphism groups (the canonical bases are explicitly given). \\

{\bf The set $\CA (m, \G )$ is an affine algebraic variety.} It is proven that each set  $\CA (m, \G )$ is an affine algebraic variety. Lemma \ref{a21Jun20}, Theorem \ref{XX12May20} and  Theorem \ref{XY12May20}
give explicit generators and defining relations for the algebra  $\OO (\CA (m, \G ))$ of regular functions on $\CA (m, \G )$
in the cases (I), (II) and (III), respectively. In the case (III), two different sets of defining relations are given (Theorem \ref{XX12May20} and  Theorem \ref{XY12May20}). Theorem \ref{XY12May20} gives a lower bound for the dimension of the algebraic variety $\CA (m , \G )$.\\

 The algebraic torus $\mT =\{ t_\l \, | \, \l \in K^\times \}$ acts on the set $\CA (m)$ by the rule:
$$\mT \times \CA (m)\ra \CA (m), \;\; (t_\l , A)\mapsto t_\l (A),$$
where $t_\l (x) = \l x$. The subsets $\CA (m, \G )$ of $\CA (m)$, where $\G \in \mS (m)$, are $\mT$-stable (that is $\mT \CA (m, \G ) =\CA (m , \G )$). By Theorem \ref{3May20}.(1), the set of $\mT$-orbits of $\CA (m)$,
$$
\mA (m):= \CA (m)/\mT = \{ \mT A\, | \, A\in \CA (m)\},
$$
 is the set of isomorphism classes of algebras  $\CA (m)$. By (\ref{grbA2}), 
$$
\mA (m)=\coprod_{\G \in \mS (m)}\mA (m , \G )\;\; {\rm where}\;\; \mA (m,\G ) := \CA (m, \G )/\mT
$$
is the set of isomorphism classes of algebras in $\CA (m, \G )$, see Section \ref{MODULISP} for details. The sets $\mA (m , \G )$ are explicitly described by 
Theorem \ref{3May20}.(2) and Proposition \ref{a20Jul20}. 

For each automorphism group $G$, Corollary \ref{a12May20} gives an explicit description  of the isomorphism classes of algebras in $\CA (m, \G )$ with automorphism group $G$.  \\

{\bf Connections with the problem of classification of singularities of holomorphic maps.} 
$K[[x]]$ is the  algebra of formal power series in the variable  $x$ over $K$; $m\geq 2$ is a natural number; $\hCA (m)$ is the set of all $K$-subalgebras $B$ of $K[[x]]$ such that $x^mK[[x]]\subset B$ but $x^{m-1}\not\in B$ i.e. the ideal 
 $x^mK[[x]]$ of $B$ is the largest ideal of $B$ which is also an ideal of $K[[x]]$.
 For each $\G \in \mS (m)$, let $\hCA (m, \G ) :=\{ B\in \hCA (m)\, | \, \G_B=\G\}$ where $\G_B$ is defined in the same fashion as in the affine case. Then 
$$\hCA (m)=\coprod_{\G\in \mS (m)}\hCA (m,\G).$$
  
The `analytic'  group of automorphisms $G^{an}=\{ x\mapsto \sum_{i\geq 1} \l_ix^i\}$ of   the algebra $K[[x]]$ respects the sets $\hCA (m, \G)$ (where $\l_i\in K$ and $\l_1\neq 0$). Classification of the orbits $\hCA (m, \G )/G^{an}$ is an old open classical problem (see, for example, \cite{Arnold-1998}), the problem of classification of curve singularities (a curve singularity  is the germ of a holomorphic map of the complex line into   complex space at a singular point).

\begin{proposition}\label{A18Jul20}
  The map $\CA (m, \G ) \ra \hCA (m, \G)$, $A\mapsto A+x^mK[[x]]$ is a bijection with inverse $B\mapsto B/x^mK[[x]]+x^mK[x]$ where $B/x^mK[[x]]\subseteq \oplus_{i=0}^{m-1}Kx^i\subseteq K[x]$. 
\end{proposition}

{\it Proof}. Straightforward. $\Box $\\

The algebra  $A+x^mK[[x]]$ is the $(x)$-adic completion of the algebra $A$. So, each algebra $B$ in $\hCA (m , \G )$ is the completion of the  algebra  $B/x^mK[[x]]+x^mK[x]$ in $\CA (m , \G )$ which is called the {\em affine partner} of the algebra $B$.
 So, the class of algebras $\hCA (m, \G )$ is completely described by the class $\CA (m ,\G )$.  
 Clearly, if two algebras in $\CA (m , \G)$ are isomorphic then so are their completions, i.e. the map  $\CA (m, \G ) \ra \hCA (m, \G)$, $A\mapsto A+x^mK[[x]]$ is isomorphism-invariant.  The classification (up to isomorphism) of algebras in $\CA (m, \G)$  is finer than the classification of topological algebras  in $\hCA (n, \G)$ up to the action of the group $G^{an}$. Furthermore,  
 
\begin{equation}\label{AmAG}
\hCA (m,\G )/G^{an}\simeq \mA (m, \G)/G^{an}_1:=\overline{\mA} (m, \G )/G^{an}_1(m)
\end{equation}
 where $G^{an}_1$ is a subgroup of $G^{an}$ that contains all automorphisms with $\l_1=1$, $\overline{\mA} (m, \G ):=\{ \bA \, |\, A\in \CA (m, \G)\}$ and      $G^{an}_1(m)$ is the subgroup of  automorphism of the algebra $K[x]/(x^m)$ of the type $x\mapsto x+\sum_{i=2}^{m-1}\l_ix^i$ where $\l_i\in K$.


\section{The canonical basis, generators and defining relations of algebras $\bA$ and $A$ where $A\in \CA (m, \G )$}\label{CANBAS}

{\em In this section, the field $K$ is an arbitrary field, i.e. not necessariuly algebraically closed}.

The aim of the section is to show that each algebra $\bA =A/x^mK[x]$, where $A\in \CA (m)$, admits a unique $K$-basis which is called the canonical basis (Proposition \ref{A4May20}) that has remarkable properties as we will see later in the paper. Existence 
 of this basis basis is a key fact in finding the automorphism group of the algebra $A$ and an isomorphism criterion for algebras $A$. Using the canonical basis of $\bA$, several explicit sets of generators and definig relations are found for the algebra $\bA$ (Corollary \ref{a5May20}, Theorem \ref{6May20} and Theorem \ref{A6May20}) and for the algebra $A$ (Theorem  \ref{Aa5May20}, Theorem \ref{X6May20} and Theorem \ref{XA6May20})

Given an algebra $A\in \CA (m)$, i.e. $K[x]\supset A\supseteq x^mK[x]$ and $x^{m-1}\not\in A$, $m\geq 2$, the polynomial algebra $K[x]$ is the normal closure of the algebra $A$ and $(x^m)=x^m K[x]$ is the largest ideal of $A$ which is also an ideal of $Kx]$. It is called the {\em conductor} of $A$ and is the sum of all ideals of $K[x]$ that are contained in $A$. Then 
$$\bA :=A/x^mK[x]\subset F:=K[x]/(x^m)$$
is a proper inclusion of finite dimensional, local, commutative $K$-algebras. In particular, their maximal ideals $A\cap (x)$ and $(x)$ are their radicals $\rad (\bA )$ and $\rad (F)$, respectively. \\

{\bf The induced radical filtration on $\bA$ and the semigroup $\G_A\cup \{ \infty\}$ of $A$.} The {\em radical filtration} (the $(x)$-{\em adic filtration}) on $F$, 
$$ F\supset (x)\supset (x)^2\supset \cdots \supset (x)^{m-1}\supset (x)^m=0,$$ induces the {\em induced radical filtration} on the algebra $\bA$, 
$$\bA\supseteq \bA_{\geq 1}\supseteq \bA_{\geq 2}\supseteq\cdots \supseteq \bA_{\geq m-1}\supseteq \bA_{\geq m}=0,$$ where $\bA_{\geq i}
:=\bA\cap (x)^i$. The asociated graded algebra of $F$, ${\rm gr}(F):=\bigoplus_{i\geq 0}(x)^i/(x)^{i+1}$ is isomorphic to the algebra $F$. 

Let $\N_+:=\{ 1, 2 , \ldots \}$. Then $(\N_+, +)$ is an abelian semigroup and the set $m+\N$ is its ideal ($\N_++(m+\N ) \subseteq m+\N$). Let $\N_{+,m}:=\N_+/(m+\N )$, a factor semigroup. Then $\N_{+,m}=\{1,2,\ldots , m-1, \infty\}$ where the addition in $\N_{+,m}$ is given by the rule
$$i+j=\begin{cases}
i+j& \text{if }i+j<m,\\
\infty& \text{if }i+j\geq m,\\
\end{cases}$$
and $i+\infty = \infty$ for all elements $i\in \N_{+,m}$. We denote by $\Sub (\N_{+.m})$ the set of all subsemigroups of $\N_{+,m}$. The algebra ${\rm gr}
(F)\simeq F=\bigoplus_{i\in \N_m}Kx^i$ is an $\N_m$-graded algebra where 
 $\N_m:=\N /(m+\N )=\{0,1,\ldots m-1, \infty\}$ is a commutative monoid, $x^\infty :=0$ and $Kx^\infty =0$. Clearly, $\N_{+,m}$ is a subsemigroup of $\N_m$. The associated graded algebra ${\rm gr}(\bA )=\bigoplus_{i=0}^{m-1}\bA_{\geq i}/\bA_{\geq i+1}$ is a homogeneous subalgebra of ${\rm gr}(F)=F$. In particular, 
 
\begin{equation}\label{grbA}
{\rm gr}(\bA )=K\oplus \bigoplus_{\g \in \G_A}Kx^\g
\end{equation}
where $\G_A:=\{ \g \, | \, x^\g \in {\rm gr}(\bA ), 1\leq \g \leq m-1\}$. Clearly, 
\begin{equation}\label{grbA1}
\G_A+\G_A\subseteq \G_A\coprod [m, \infty )\;\; {\rm and}\;\; m-1\not\in \G_A.
\end{equation}
In particular, the set $\G_A\cup\{ \infty\}$ is a subsemigroup of $\N_{+,m}$. Notice that $\G_A=\emptyset$ iff $A=K+x^mK[x]$. We will see that the number $m$ and $\G_A$ are isomorphism-invariants of the algebra $A$ (algebras $A\in \CA (m)$ and $A'\in \CA (m')$ are isomorphic then $m=m'$ and $\G_A =\G_{A'}$, Theorem \ref{3May20}.(2)).

Let $\mS (m)$ be the set that contains the empty set and all non-empty subsets $\G$ of the set $\{ 2,\ldots , m-2\}$ such that $\G +\G \subseteq \G \cup [m, \infty )$. Notice that, by definition, $m-1\not\in \G$; $\mS (m)\neq \emptyset$ iff $m\geq 4$ (for all $m\geq 4$, $\G =\{ m-2\}\in \mS (m)$); and the map from the set $\mS (m)$ to the set of all subsemigroups of $\N_{+, m}$ which is given by the rule $\G \mapsto \G \cup \{ \infty\}$, is an {\em injection}.

For each $\G \in \mS (m)$, let $\CA (m, \G ) :=\{ A\in \CA (m)\, | \, \G_A=\G\}$ where $\G_A$ is the set of all natural numbers $i$ such that $2\leq i<m-1$ and there is an element $a_i=x^i+\sum_{j>i}\l_{ij}x^j\in A$ for some $\l_{ij}\in K$. Clearly, $\G_A\in \mS (m)$.  Then 
\begin{equation}\label{grbA2}
\CA (m)=\coprod_{\G\in \mS (m)}\CA (m,\G).
\end{equation}
We will see that each set $\CA (m, \G )$ is an affine algebraic variety (Lemma \ref{a21Jun20}, Theorem \ref{XX12May20} and  Theorem \ref{XY12May20}). \\

{\bf The canonical basis of the algebra $\bA$ where $A\in \CA (m, \G )$.} Let $C\G =\{ 2, \ldots , m-1\} \backslash \G$. For each $i\in \{ 2, \ldots , m-1\}$, let $\G (i):=\{ \g \in \G \, | \, \g >i\}$ and $C\G (i) :=\{ \d \, | \, \d\not\in \G, i<\d \leq m-1\}$. 

 \begin{proposition}\label{A4May20}
For each algebra $A\in \CA (m, \G \}$, there is a unique basis $\{ 1, f_\g, \, | \, \g \in \G \}$ of the algebra $\bA$ that satisfies the condition that $f_\g =x^\g +\sum_{\d \in C\G (\g )}\l_{\g \d}x^\d$ where $\l_{\g \d}\in K$ and $ C\G (\g ) =\{ \d \, | \, \d\not\in \G, \g <\d \leq m-1\}$. 
\end{proposition}

{\it Proof}. {\em Existence:} By (\ref{grbA}), we can find a $K$-basis $\{ 1, g_\g \, | \, \g \in \G\}$ where $g_\g =x^\g +\sum_{\g <\d \leq m-1}Kx^\d$ and $\l_{\g \d }\in K$. Suppose that $\G =\{ \g_1, \ldots , \g_t\}$ where $2\leq \g_1< \cdots  <\g_t\leq m-2$. The element  $f_{\g_t}:=g_{\g_t}$ satisfies the condition of the proposition. Now, we use the downward  induction on $i$ starting at $i=t$. To show existence we  suppose that $i<t$ and we have  found already the elements $f_{\g_{i+1}}, \ldots , f_{\g_t}$ that satisfy the condition  given in the proposition. The element $g_{\g_i}$ can be written as follows
$$ g_{\g_i}=x^{\g_i}+\sum_{\d \in C\G (\g )} \l_{\g_i\d}x^\d +\sum_{\g'\in \G (\g )}\l_{\g_i \g'}x^{\g'}.$$
Then 
$$f_{\g_i}:=g_{\g_i}-\sum_{\g'\in \G (\g )}\l_{\g_i\g'}f_{\g'}$$ satisfies the condition of the proposition. Now,  existence follows by induction on $i$.
 
 (ii) {\em Uniqueness:} If $\G =\emptyset$ theoren there is nothing to prove. Suppose that $\G \neq \emptyset$ and $ \{ 1, f_\g'\, | \, \g \in \G\}$ is another $K$-basis as in the proposition. Then 
 $$ f_\g -f_\g'\in \bA\cap \Bigg(\sum_{\d \in C\G (\g )}Kx^\d \Bigg)=0$$
 for all elements $\g \in \G$, and we are done.  
$\Box $ \\

{\it Definition.} The $K$-basis $\{ 1, f_\g \, | \, \g \in \G\}$ in Proposition \ref{A4May20} is called  {\bf the canonical $K$-basis} of the algebra $\bA$. Any $K$-basis of the type $\{ 1, \l_\g f_\g \, | \, \g \in \G\}$, where $\l_\g\in K^\times$, is called 
 {\bf a canonical $K$-basis} of the algebra $\bA$.\\
 
 {\bf The structure constants of the algebra $\bA$ with respect to  the canonical basis.} For an algebra $A\in \CA (m, \G )$, let $\{1, f_\g \, | \, \g \in \G \}$ be the  canonical basis of the algebra $\bA =A/x^mK[x]$ where 
\begin{equation}\label{fgg}
f_\g = x^\g +\sum_{\d \in C\G (\g )}\l_{\g \d}x^\d\;\; {\rm and}\;\; \l_{\g\d}\in K.
\end{equation}
The total number of the parameters $\{ \l_{\g\d}\}$ is equal to 
\begin{equation}\label{fgg1}
\prod_{\g \in\G } |C\G (\g )|.
\end{equation}
In general, they are not algebraically independent. For each subset $S$ of the set $\{ 2, \ldots , m-1\}$, let us consider the {\em characteristic function} of the set $S$,
$$\chi (\cdot , S):\{ 2, \ldots , m-1\}\ra \{0,1\}, \;\; i\mapsto \chi (i,S)=
\begin{cases}
1& \text{if } i\in S,\\
0& \text{if }i\not\in S.\\
\end{cases}$$
Then, by direct computation, for all elements $\g , \g'\in \G$, 

\begin{equation}\label{fgg2}
f_\g f_{\g'}=
\begin{cases}
x^{\g +\g'}+\sum_{\xi\in C\G (\g +\g')}\l_{\g , \g' ; \xi}x^\xi+\sum_{\rho\in \G (\g +\g')}\mu_{\g , \g';\rho}x^\rho & \text{if }\g +\g'<m-1,\\
0& \text{if }\g +\g'\geq m-1.\\
\end{cases}
\end{equation}
where 
\begin{equation}\label{fgg2l}
\l_{\g , \g' ; \xi}=\l_{\g, \xi-\g'}\chi (\xi-\g', C\G (\g ))+\l_{\g', \xi-\g}\chi(\xi-\g , C\G (\g'))+\sum_{\{ \d+\d'=\xi\, | \, \d\in C\G (\g ), \d'\in C\G (\g')\} }\l_{\g\d}\l_{\g'\d'},
\end{equation}

\begin{equation}\label{fgg2m}
\mu_{\g , \g' ; \rho}=\l_{\g, \rho-\g'}\chi (\rho-\g', C\G (\g ))+\l_{\g', \rho-\g}\chi(\rho-\g , C\G (\g'))+\sum_{\{ \d+\d'=\rho\, | \, \d\in C\G (\g ), \d'\in C\G (\g')\} }\l_{\g\d}\l_{\g'\d'}.
\end{equation}
For elements $i, j\in \{2, \ldots , m-1\}$ such that $i<j$, let $\G (i,j):=\{k\in \G\, | \, i<k<j\}$ and $C\G (i,j):=\{k\in C\G\, | \, i<k<j\}$. 

For each algebra $A\in \CA (m, \G )$, Proposition \ref{A5May20}.(1) and (\ref{fgg2m}) determine explicitly the structure constants of the algebra $\bA$ with respect  to its canonical basis. 

\begin{proposition}\label{A5May20}
Given an algebra $A\in \CA (m, \G )$. Let $\{ 1, f_\g \, | \,\ \g \in \G \}$ be the canonical basis of the algebra $\bA$, see (\ref{fgg}). Then for all elements $\g , \g'\in \G$ such that $\g+\g'<m-1$, 
\begin{enumerate}
\item $f_\g f_{\g'}=f_{\g +\g'}+\sum_{\rho \in \G (\g +\g')}\mu_{\g , \g';\rho} f_\rho$, and 
\item $f_{\g+\g'}=x^{\g +\g'}+\sum_{\xi \in C\G (\g +\g')}\l_{\g, \g';\xi}x^\xi -\sum_{\d \in C\G (\g +\g')}\Bigg( \sum_{\rho \in \G (\g+\g', \d )}\mu_{\g , \g';\rho}\l_{\rho \d}\Bigg) x^\d$. 
\item For all elements $\d\in C\G (\g +\g')$, $\l_{\g +\g', \d}=\l_{\g , \g';\d}-\sum_{\rho \in \G (\g+\g', \d )}\mu_{\g , \g';\rho}\l_{\rho \d}$.
\end{enumerate}
\end{proposition}

{\it Proof}.  Let $D$ be the double sum in statement 2. By (\ref{fgg}) and (\ref{fgg2}), 

\begin{eqnarray*}
 f_\g f_{\g'}&=&x^{\g +\g'}+\sum_{\xi \in C\G (\g +\g')}\l_{\g , \g';\xi} x^\xi -\sum_{\rho \in \G (\g +\g')}\sum_{\d\in C\G (\rho )}\mu_{\g , \g';\rho} \l_{\rho\d}x^\d +\sum_{\rho \in \G (\g +\g')}\mu_{\g , \g';\rho} f_\rho\\
 &=& x^{\g +\g'}+\sum_{\xi \in C\G (\g +\g')}\l_{\g , \g';\xi} x^\xi -D+ \sum_{\rho \in \G (\g +\g')}\mu_{\g , \g';\rho} f_\rho\\
 &=&f_{\g+\g'}+\sum_{\rho \in \G (\g +\g')}\mu_{\g , \g';\rho} f_\rho +\Bigg\{ -\sum_{\d\in C\G (\g +\g')}\l_{\g+\g', \d}x^\d +\sum_{\xi \in C\G (\g+\g')}\l_{\g,\g';\xi}x^\xi-D\Bigg\}.
\end{eqnarray*}
Since $\bA \ni f_\g f_{\g'}-f_{\g+\g'}-\sum_{\rho \in \G (\g +\g')}\mu_{\g , \g';\rho} f_\rho =\{ \cdots \}$, we must have $\{ \cdots \}=0$ (as three sums in $\{ \cdots \}$ are over $C\G (\g +\g')$ and the LHS of the equality is an element of the algebra $\bA$ and if it is nonzero then the LAH $= \l x^{\g''} +\cdots $ for some element $\l \in K^\times$ and $\g''\in \G$, a contradiction), and statements 1--3 follow. $\Box $\\

{\bf Generators and defining relations of the algebra $\bA$.} Corollary \ref{a5May20} defines the commutative algebra $\bA$ via generators and defining relations.

\begin{corollary}\label{a5May20}
Given an algbera $A\in \CA (m, \G )$, we keep the notation of Proposition \ref{A5May20}. Then the commutative algebra $\bA$ is generated by the elements $\{ f_\g\}_{\g \in \G}$ subject to the defining relations:
\begin{eqnarray*}
 f_\g f_{\g'}&=& f_{\g+\g'}+\sum_{\rho \in \G (\g +\g')}\mu_{\g, \g';\rho}f_\rho\;\; {\rm for\; all}\;\; \g , \g'\in \G\;\; {\rm such \; that}\;\; \g+\g'<m-1,\\
 f_\g f_{\g'}&=& 0\;\; {\rm for\; all}\;\; \g , \g'\in \G\;\; {\rm such \; that}\;\; \g+\g'\geq m-1.
\end{eqnarray*}
\end{corollary}

{\it Proof}. The corollary follows at once from Proposition \ref{A5May20}.(1) and (\ref{fgg2}). $\Box $ \\

Each non-empty set $\G \in \mS (m)$ is a  disjoint union $$\G = \ind (\G ) \coprod \dec (\G )$$ where $\ind (\G )$ is the set of {\em indecomposable} elements of $\G$ and $\dec (\G )$ is the set of {\em decomposable} elements of $\G$. An element $\g \in \G$ is called {\em decomposable} if $\g =\g_1+\g_2$ for some elements $\g_1, \g_2\in \G$. Fix $\G \in \mS (m)$ and let 
\begin{equation}\label{indG}
\ind (\G ) =\{ \nu_1, \ldots , \nu_s\}\;\; {\rm where}\;\; 2\leq \nu_1< \cdots < \nu_s\leq m-2. 
\end{equation}
Let $\Q$ be the field of rational numbers. For the $s$-dimensional vector space $\Q^s=\bigoplus_{i=1}^s\Q e_i$, let $ab:= \sum_{i=1}^sa_ib_i$ the {\em scalar product} on it where $a=(a_1, \ldots , a_s), b=(b_1, \ldots  , b_s)\in \Q^s$. The set of elements $e_1, \ldots , e_s$ is the canonical $\Q$-basis of $\Q^s$. Define 
\begin{equation}\label{nuQ}
\nu_\G: =\nu :=(\nu_1, \ldots , \nu_s) \in \N_+^s\subset  \Q^s.
\end{equation}
For each element $\g \in \dec (\G )$, let 
$$\Rel (\g ) :=\{ a\in \N^s\, | \, a\nu = \g \}.$$
Clearly $|\Rel (\g )|<\infty$. For each element $\g \in \dec (\G )$, we fix an element
\begin{equation}\label{nuQ2}
a(\g ) \in \Rel (\g ).
\end{equation}
For each $\nu_i\in \ind (\G )$, we set $a(\nu_1) = \nu_i$.  In particular, $a(\g ) \nu = \g$ for all $\g \in \G$. The choice of $a(\g )$ is arbitrary but fixed. By definition, $a(\nu_i) := e_i$ for all $i=1, \ldots , s$.  Let 
\begin{equation}\label{nuQ1}
\dec (\G )_{\geq 2}:=\{ \g \in \dec (\G )\, | \, |\Rel (\g )|\geq 2\}.
\end{equation}
Recall that $f_{\nu_i}=x^{\nu_i}+\sum_{j\in C\G (\nu_i)}\l_{\nu_i,j}x^j$ for $i=1, \ldots , s$ and   $C\G (\nu_i, k)=\{ j\, | \, \nu_i<j<k, j\not\in \G \}$. 
For each natural number $k$ such that $2\leq k\leq m-1$, let 
\begin{equation}\label{Lgdef}
\L (k) :=\{ \l_{\nu_i,j}\, | \, i=1, \ldots , s; \nu_i<k; j\in C\G (\nu_i, k)\}.
\end{equation}

Given an algebra $A\in \CA (m, \G )$.  Let $\{ 1, f_\g \, | \, \g \in \G \}$ be the canonical basis of $A$. For each $a=(a_i)\in \N^s$, let $f^a:=\prod_{i=1}^sf_{\nu_i}^{a_i}$. In particular, for each $\g \in \dec (\G )$ and $\nu_i\in \ind (\G )$, $f^{a(\g )}=\prod_{i=1}^sf_{\nu_i}^{a(\g )_i}$ and $f^{a(\nu_i)}=f_{\nu_i}$. 


\begin{lemma}\label{a16May20}
For all nonzero elements $a\in \N^s$, 
\begin{eqnarray*}
 f^a&=&f_{a\nu}+\sum_{\g'\in \G (a\nu )}c_{\g'}(f^a)f_{\g'}, \\
 f_{a\nu}&=&x^{a\nu}+\sum_{a\nu<\d <m}c_\d (f^a)x^\d -\sum_{\g'\in \G (a\nu )}c_{\g'} (f^a)f_{\g'},
\end{eqnarray*}
 where $c_{\g'}(f^a)$ is the coefficient of $x^{\g'}$ of the polynomial $f^a=x^{a\nu }+\sum_{j=a\nu +1}^{m-1}c_j(f^a)x^j$. The coefficient $c_{\g'}(f^a)$ is an explicit polynomial in the variables $\L (\g')$, see (\ref{Lgdef}).
\end{lemma}

{\it Proof}. Recall that $f_{\nu_i}=x^{\nu_i}+\sum_{j\in C\G (\nu_i)}\l_{\nu_i,j}x^j$ for $i=1, \ldots , s$. By multiplying out, we get the equality  $$f^a=x^{a\nu } + \sum_{j=a\nu +1}^{m-1} c_j(f^a) x^j$$ where each $c_j(f^a)$ is an explicit polynomial in the variables $\L (j)$. For each element $\g \in \G$, let $t_\g :=f_\g -x^\g$, the tail of the polynomial $f_\g$. Then 
\begin{eqnarray*}
 f^a&=&x^{a\nu} +\sum_{\d \in C\G (a\nu )}c_\d (f^a)x^\d -\sum_{\g'\in \G (a\nu )}c_{\g'} (f^a)t_{\g'}+\sum_{\g'\in \G (a\nu )}c_{\g'}(f^a)f_{\g'}\\
 &=& f_{a\nu}+\sum_{\g'\in \G (a\nu )}c_{\g'}(f^a)f_{\g'} +\D\\ 
\end{eqnarray*}
where 
$ \D = -t_{a\nu}+\sum_{\d \in C\G (a\nu )}c_\d (f^a)x^\d -\sum_{\g'\in \G (a\nu )}c_{\g'} (f^a)t_{\g'}\in \sum_{\d\in C\G (a\nu)}Kx^\d$. Then $$\D=0$$ since 
$\D= f^a-f_{a\nu}-\sum_{\g'\in \G (a\nu )}c_{\g'}(f^a)f_{\g'}\in \bA$ (if $\D \neq 0$ then on the one hand $\D \in K^\times x^{\d}+\cdots $ for some $\d\in \G (a\nu )$ since $\D \in \bA$, on the other hand $\d\in C\G (a\nu )$, a contradiction; the three dots denote higher degree terms), and the first equality of the lemma follows. 

The second equality of  the lemma follows from the first and the third equalities of the lemma (the second equality of the lemma is equivalent to the equality $\D =0$). $\Box$\\

{\bf Elements of the canonical basis as polynomial in $f_{\nu_1}, \ldots , f_{\nu_s}$.}
Lemma \ref{a6May20}.(2) represents  each element $f_\g$ of the canonical basis of the algebra $A\in \CA (m , \G )$ as an explicit polynomial in $f_{\nu_1}, \ldots , f_{\nu_s}$. For each element $\g \in \dec (\G )_{\geq 2}$ and element $b\in \Rel (\g ) \backslash \{ a(\g )\}$, Lemma \ref{a6May20}.(3) represents  the element $f^b$ as an explicit linear combination of the elements $\{ f^{a(\g )}\, | \, \g \in \G\}$. 

\begin{lemma}\label{a6May20}
Let $A\in \CA (m , \G )$ and $\ind (\G ) = \{ \nu_1, \ldots , \nu_s\}$ (see above). Then 
\begin{enumerate}
\item The set $\{ 1, f_{\nu_1}, \ldots , f_{\nu_s}, f^{a(\g )}\, | \, \g \in \dec (\G )\}$ is a $K$-basis of the algebra $\bA$. 
\item For each element $\g \in \G$, $f_\g = f^{a(\g )}+\sum_{\g'\in \G (\g )}\eta_{\g\g'}f^{a(\g')}$ for unique elements $\eta_{\g\g'}\in K$. Each element $\eta_{\g\g'}$ is an explicit polynomial in $\L (\g')$, see (\ref{etagg}) (if $\g = \nu_i$ then  $f_{\nu_i}=f^{a(\nu_i)}=f_{\nu_i}$ is a tautology).
\item For each element $\g \in \dec (\G )_{\geq 2}$ and each element $b\in \Rel (\g ) \backslash \{ a(\g )\}$, $$f^b=f^{a(\g )}+\sum_{\g'\in \G (\g )}\th_{\g , \g ';b}f^{a(\g')}$$ for unique elements $\th_{\g, \g';b}\in K$. Each element $\th_{\g , \g'; b}$ is an explicit polynomial in $\L (\g')$, see also (\ref{thggb}).  
\end{enumerate}
\end{lemma}

{\it Proof}. 1. Statement 1 follows from the fact that $f^{a(\g )}\in f_\g +\sum_{\g'\in \G (\g )}Kf_{\g'}$ where $\{ 1, f_\g \, | \, \g \in \G\}$ is the canonical basis of the algebra $\bA$. 

2. Statement 2 follows from statement 1 and Lemma \ref{a16May20}.

3. Statement 3 follows from statement 2 and Lemma \ref{a16May20}. $\Box $\\

{\bf Explicit expressions for the elements $\eta_{\g\g'}$ and $ \th_{\g, \g';b}.$}
 We keep the notation of Lemma \ref{a6May20}.   In order to write explicit defining relations  of the algebra $\bA$  (Theorem \ref{X6May20}) and the algebra $A$ (Theorem \ref{XA6May20}) we need an explicit expression for the elements $ \th_{\g, \g';b}$, see (\ref{thggb}). In order to obtain it we use Lemma \ref{a16May20} and Lemma \ref{a6May20}. There are two $K$-bases for the algebra $\bA$: the standard basis $\{ 1, f_\g\, | \, \g \in \G \}$ where the order of elements are increasing (if $\G =\{ \g_1<\cdots <\g_t\}$ then $\{ 1<f_{\g_1}<\cdots <f_{\g_t}\}$) and the $K$-basis $\{ 1, f^{a(\g )}\, | \, \g \in \G\}$ where the order of elements is also increasing, i.e. 
 $\{ 1<f^{a(\g_1)}<\cdots <f^{a(\g_t)}\}$.
 Recall that $a(\nu_i)=e_i$ for all $i=1, \ldots , s$. By Lemma \ref{a16May20}, 
 \begin{eqnarray*}
 f^{a(\g)}&=&f_\g+\sum_{\g'\in \G (\g)}c_{\g'}(f^{a(\g )})f_{\g'}, \;\;\g \in \dec (\G),\\
  f^{a(\nu_i)}&=& f_{\nu_i}, \;\; i=1, \ldots , s.
\end{eqnarray*}
 It follows from the equalities above that the change-of-basis matrix $C_{\bA}$ (from the standard basis to the second one) is a {\em lower unitriangular matrix} (i.e. a lower triangular matrix with identities on the diagonal). The elements of the matrix $C_{\bA}$ are the coefficients in the equalities above. So, the matrix $C_{\bA}$ is a sum $C_{\bA}=1+\gn_{\bA}$ of the identity matrix 1 and a strictly lower triangular matrix $\gn_{\bA}$. Hence, the inverse matrix of $C_{\bA}$, 
\begin{equation}\label{etagg}
 C_{\bA}^{-1}=\sum_{i=0}^{\dim (\bA )} (-\gn_{\bA})^i,
\end{equation} 
  is also a lower unitriangular matrix elements of which are explicit expressions, i.e. {\em the elements $\eta_{\g\g'}$ (in Lemma \ref{a6May20}.(2)) are explicit expressions.}  
 Now, when we substitute the expression for $f^{a(\g )}$ from Lemma \ref{a16May20} into the sum in Lemma \ref{a6May20}.(2) to obtain the result: For each element $\g \in \dec (\G )_{\geq 2}$ and element $b\in \Rel (\g ) \backslash \{ a(\g )\}$,
 \begin{eqnarray*}
 f^b&=& f_\g +\sum_{\g'\in\G (\g )}c_{\g'}(f^b)f_{\g'}\\
 &=&f^{a(\g )}+\sum_{\g'\in \G (\g )}\eta_{\g \g '}f^{a(\g')}+\sum_{\g'\in\G (\g )}c_{\g'}(f^b)\bigg(f^{a(\g' )}+\sum_{\g''\in \G (\g' )}\eta_{\g' \g ''}f^{a(\g'')}\bigg)\\
 &=& f^{a(\g )}+\sum_{\g'\in \G (\g )}\bigg( \eta_{\g\g'}+c_{\g'}(f^b)+\sum_{\d \in \G (\g , \g' )}c_\d (f^b)\eta_{ \d\g'}
 \bigg)f^{a(\g')}.
\end{eqnarray*}
In order to get the last equality we changed the order of summation in the double sum and replaced $(\g', \g'')$ by $(\d, \g')$. 
 Therefore, for each element $\g \in \dec (\G )_{\geq 2}$ and element $b\in \Rel (\g ) \backslash \{ a(\g )\}$,
\begin{equation}\label{thggb}
\th_{\g , \g'; b}=\eta_{\g\g'}+c_{\g'}(f^b)+\sum_{\d \in \G (\g , \g' )}c_\d (f^b)\eta_{\d\g'}.
\end{equation}

By Lemma \ref{a6May20}.(3), the elements $\th_{\g , \g'; b}$ can be also found recursively by using the equality
\begin{equation}\label{thggb1}
\th_{\g , \g'; b}=c_{\g'}\bigg( f^b-f^{a(\g )}-\sum_{\g''\in \G (\g , \g' )}\th_{\g , \g '';b}f^{a(\g'')}\bigg)
\end{equation}
where the RHS is the coefficient of $x^{\g'}$ of the polynomial in the brackets.\\

Consider the monoid $(\N^s,+)$ and its subsemigroup 
$$\N^s(m-1,\G ):=\{ c\in \N^s\, | \, c\nu\geq m-1\}=\{ c\in \N^s\, | \, c\nu > m-1\}.$$ 
The subsemigroup $\N^s(m-1,\G )$ is also an ideal of the monoid $(\N^s,+)$, that is $\N^s +\N^s(m-1,\G )\subseteq \N^s(m-1,\G )$. The set 
$$\ind \, \N^s(m-1,\G ):=\N^s(m-1,\G )\backslash \Big(\N^s(m-1,\G )+(\N^s\backslash \{ 0\} )\Big)$$ of minimal generators of the ideal $\N^s(m-1,\G )$ is a {\em finite set}. 

{\it Proof.} Let $\L$ be the (homogeneous) subalgebra of $K[x]$ that is generated by the elements $x^{\nu_1}, \ldots , x^{\nu_s}$ (where $\ind (\G ) = \{ \nu_1, \ldots , \nu_s\}$). Then the ideal $ \ga =\L\cap x^{m-1}K[x]$ of the Noetherian algebra $\L$ is a finitely generated and the finite set $\{x^{\d \nu}\, | \, \d\in \ind \, \N^s(m-1, \G)\}$ is its {\em minimal} set of generators (as a $\L$-module). In particular, the factor algebra
\begin{equation}\label{thggb2}
K[x^{\nu_1}, \ldots , x^{\nu_s}]/(x^{\d \nu}\, | \, \d\in \ind \, \N^s(m-1, \G))=K\oplus\bigoplus_{\g \in \G}Kx^\g 
\end{equation}
is finite dimensional. Its dimension is $1+|\G |$. $\Box$ \\

Lemma \ref{a7May20} is a characterization of algebras $A\in \CA (m)$ with $|\ind (\G_A )|=1$. For a rational number $q\in\Q$, we denote by $[q]$  the integer part of $q$, i.e. the unique integer $[q]$ such that $[q]\leq q<[q]+1$. 

\begin{lemma}\label{a7May20}
Given an algebra $A\in \CA (m, \G )$. Suppose that $\ind (\G ) =\{ \nu_1\}$. Then $\bA \simeq K[f_{\nu_1}]/(f_{\nu_1}^{[\frac{m-1}{\nu_1}]+1})$ and vice versa, i.e. if $\bA \simeq K[f]/(f^{[\frac{m-1}{d}]+1})$ for some algebra $A\in \CA (m)$ where $f=x^d+\l_{d-1}x^{d-1}+\cdots +\l_{m-1}x^{m-1}\in K[x]$ then $\ind (\G ) =\{ d\}$, $d\geq 2$  and $d\nmid m-1$. 
\end{lemma}

{\it Proof}. $(\Rightarrow )$ Since $\ind (\G ) = \{ \nu_1\}$, $\G =\{ i\nu_1\, | \, 1\leq i\leq [\frac{m-1}{\nu_1}]\}$ and $\bA =K\oplus\bigoplus_{i=1}^{[\frac{m-1}{\nu_1}]}Kf_{\nu_1}^i$ with $f_{\nu_1}^{[\frac{m-1}{\nu_1}]+1}=0$, and the implication follows. 

$(\Leftarrow )$ Straightforward. $\Box $\\

In view of Lemma \ref{a7May20}, we will assume that $s=|\ind (\G )|\geq 2$. Theorem \ref{6May20} shows that the set $\{ f_{\nu_1}, \ldots , f_{\nu_s}\}$ is a minimal set of generators of the algebra $\bA$, it also presents defining relations of the algebra $\bA$ that the set of generators  satisfies. In general, some of the defining relations that are given in the first equation of Theorem \ref{6May20} are redundant. An irredundant set is given in Theorem \ref{A6May20}.

\begin{theorem}\label{6May20}
Given an algebra $A\in \CA (m , \G )$. Suppose that $s=|\ind (\G )\geq 2$ and we keep the notation as above. Then the algebra $\bA$ is generated by the the minimal  set  of generators  $\{ f_{\nu_1}, \ldots , f_{\nu_s}\}$ that satisfies the defining relations (as a commutative algebra):
\begin{eqnarray*}
 f^b&=& f^{a(\g )}+\sum_{\g'\in \G (\g )}\th_{\g , \g';b}f^{a(\g')}\;\; {\rm where  }\;\; \g \in \dec (\G )_{\geq 2}\;\; {\rm and}\;\; b\in \Rel(\g )\backslash \{ a(\g )\},\\
f^c &=&0\;\; {\rm where  }\;\; c\in \ind \, \N^s(m-1, \G ).
\end{eqnarray*}
\end{theorem}

{\it Proof}. By Lemma \ref{a6May20}.(2), the elements $ f_{\nu_1}, \ldots , f_{\nu_s}$ are generators of the algebra $\bA$. They are a minimal set of generators for the algebra $\bA$ since $\ind (\G ) = \{ \nu_1, \ldots , \nu_s\}$. Clearly, they satisfy the relations of the theorem (by Lemma \ref{a6May20}.(3) and the definition of the semigroup $\N^s(m-1, \G )$).

Let $\bA'$ be a commutative algebra that is generated by  indeterminates $f_{\nu_1}, \ldots , f_{\nu_s}$ subject to the defining relations of the theorem. Then $\bA'=K+\sum_{\g \in \G} Kf^{a(\g )}$, by the definition of the sets $\ind \, \N^s(m-1, \G )$ and  $\dec (\G)$. Hence the epimorphism  
$$ \v : \bA'\ra \bA , \;\; f_{\nu_i}\mapsto f_{\nu_i}$$ is an isomorphism since $\bA=K\oplus\bigoplus_{\g \in \G} Kf^{a(\g )}$. $\Box$ \\

 
 {\it Remarks.} 1. The relations $`f^c=0'$ are obvious ones. They imply that the  algebra that satisfies them is a finite dimensional $\G$-graded algebra. 
 
 2. The relations $`f^b=f^{a(\g )}+\cdots '$ are not obvious ones. In general, some of them are redundant. Theorem \ref{A6May20} replaces this set of defining relation by the one none element of which can be dropped. In order to prove Theorem \ref{A6May20}, we need to introduce more concepts and to prove some more results (Lemma \ref{b7May20}). 
 
 Recall that $\nu =(\nu_1, \ldots , \nu_s)$ (see (\ref{nuQ})) and $s\geq 2$. The kernel of the $\Q$-linear map 
 $$\cdot \nu : \Q^s\ra \Q, \;\; a\mapsto a\nu=\sum_{i=1}^sa_i\nu_i$$
 is equal to $\bigoplus_{i=2}^s\Q(\nu_ie_1-\nu_1e_i)$. Its $\Q$-dimension is equal to $s-1$. The rank of the $\Z$-module $\ker_{\Q^s}(\cdot \nu )\cap \Z^s$ is equal to $s-1$ since the intersection contains $\bigoplus_{i=2}^s\Z (\nu_ie_1-\nu_1e_i)$. Consider the `ideal of relations' of the semigroup $\G$, 
\begin{equation}\label{IGdef}
I_\G :=\sum_{\g \in \dec (\G )_{\geq 2}}\sum_{a,b\in \Rel (\g )}\N (a-b)=\sum_{\g \in \dec (\G )_{\geq 2}}\sum_{a,b\in \Rel (\g )}\Z (a-b)\subseteq \ker_{\Q^s}(\cdot \nu )\cap \Z^s.
\end{equation}
 The second equality follows from the equality $\N (a-b)+\N (b-a)=\Z (a-b)$. Clearly, the set $I_\G$ is a free abelian group of rank $r=r_\G$ and $0\leq r\leq s-1$. The rank $r_\G$ is equal to zero iff $\dec (\G )_{\geq 2}=\emptyset$ iff each elemenet $\g \in \G$ is a unique sum of the  elements $\nu_1, \ldots , \nu_s$ (counted with multiplicity), and Lemma \ref{b7May20} follows from Theorem \ref{6May20}. 
 
 \begin{lemma}\label{b7May20}
Given an algebra $A\in \CA (m , \G )$. Suppose that $\dec (\G )_{\geq 2}=\emptyset$. Then $$\bA = K[f_{\nu_1}, \ldots , f_{\nu_s}]/(f^c)_{c\in \ind \, \N^s (m-1, \G )}, $$
the sets $\{f_{\nu_1}, \ldots , f_{\nu_s} \}$ and $\{ f^c \, | \, c\in \ind \, \N^s (m-1, \G )\}$ are  minimal sets of generators and defining relations of the algebra $\bA$.
\end{lemma}

In view of Lemma \ref{b7May20}, we will assume that $s=|\ind (\G )|\geq 2$ and $\dec (\G )_{\geq 2}\neq \emptyset$. \\


{\bf Sets of generators and defining relations of the algebra $\bAmon$ where $\Amon \in \CA (m, \G )$ is the monomial algebra.}  An algebra $A\in \CA (m)$ is called a {\bf monomial algebra} if it admits a monomial basis. Let $\G = \G_A$. The algebra $A$ is a monomial algebra iff $$ A=K\oplus \bigoplus_{\g \in \G} Kx^\g \oplus x^mK[x]$$ iff $ \{ 1, x^\g \, | \, \g \in \G\}$ is the canonical basis of the algebra $\bA$. We denote the monomial algebra in $\CA (m, \G )$ by $\Amon = \Amon (\G)$.

If either $|\ind (\G )|=1$ or $|\ind (\G )|\geq 2$ and $\dec (\G )_{\geq 2}=\emptyset$ then the generators and defining relations of the algebra $\bAmon$ are given 
in Lemma \ref{a7May20} and  Lemma \ref{b7May20}, respectively.

Suppose that $s=|\ind (\G )|\geq 2$ and $ \dec (\G )_{\geq 2}\neq \emptyset$. Recall that $\ind (\G) =\{ \nu_1, \ldots , \nu_s\}$. The set $f_{\nu_1} =x^{\nu_1}, \ldots ,f_{\nu_s}=x^{\nu_s}$ is a minimal set of generators of the monomial algebra $\bAmon$. Let $\mM = \mM (\G)$ be the factor algebra of the (abstract) polynomial algebra $K[f_{\nu_1},  \ldots ,f_{\nu_s}]$ in $s$ indeterminates  modulo the ideal $(f^c)_{c\in \ind \, \N^s (m-1, \G )}$. The algebra 
$$\mM = K\oplus \bigoplus_{\g \in \G} \mM_\g, \;\;{\rm where}\;\;  \mM_\g := \bigoplus_{a\in \N^s,  a\nu =\g}Kf^a, $$
is a finite dimensional, local,   graded algebra where $\gn = (f_{\nu_1},  \ldots , f_{\nu_s})$ is its maximal ideal which is a homogeneous ideal. There is a graded $K$-algebra epimorphism $\mM \ra \bAmon$, $f_{\nu_i}\mapsto f_{\nu_i}=x^{\nu_i}$ for $i=1, \ldots , s$. By Theorem \ref{6May20}, the kernel $\ga =\ga_\G$ is a homogeneous ideal which is generated by the set $\{ f^b-f^{a(\g )}\, | \, \g \in \dec (\G)_{\geq 2}, b\in \Rel (\g )\backslash \{ a(\g )\} \}$. Furthermore, 
\begin{equation}\label{BGbas1}
\mM = K\oplus\bigoplus_{i=1}^s Kf_{\nu_i} \oplus\bigoplus_{\g \in \dec (\G )}f^{a(\g )}\oplus \ga_\G \;\; {\rm and}\;\; \ga_\G= \bigoplus_{\g \in \dec_{\geq 2}(\G )}\bigoplus_{b\in \Rel (\g )\backslash \{ a(\g )\} }K(f^b-f^{a(\g )}).
\end{equation}
By  Nakayama's Lemma, any $K$-basis of the vector space  $\ga / \gn \ga$, say 
\begin{equation}\label{BGbas}
\mB = \mB_\G=\{ f^{b_1}-f^{a (\mu_1)}, \ldots , f^{b_t}-f^{ a (\mu_t)}\}, \;\; \mu_1\leq \cdots \leq \mu_t, \;\; t =t_\G :=\dim_K(\ga / \gn \ga), 
\end{equation}
 together with the set $\{ f^c\, | \, c\in \ind \, \N^s (m-1, \G ) \}$ is a set of defining relations  for the algebra $\bAmon$. 
The choice of the elements $\mu_1\leq \cdots \leq \mu_t$ (counted with multiplicity) 
 is {\em unique} since the ideals $\ga$, $\gn$  and $\gn\ga$ are homogeneous. 
 
 The elements $b_1, \ldots , b_t$ can be effectively found. There is not much freedom in their choice. 
 
 Let $\CB_\G :=\{ b\, | \, b\in \Rel (\g ) , \g \in \dec_{\geq 2}(\G )\}$, 
\begin{equation}\label{BGbas2}
 \CB_\G':=\CB _\G \backslash \Bigg(\bigcup_{i=1}^s(e_i+\CB_\G )\Bigg)=\{ b_1',  \ldots , b_\tau'\} \;\; {\rm and}\;\;   \mu_1' :=b_1'\nu ,  \ldots , \mu_\tau':=b_\tau'\nu 
\end{equation}
where $e_1, \ldots , e_s$ is the standard $\Q$-basis of the $s$-dimensional $\Q$-vector space $\Q^s=\bigoplus_{i=1}^s\Q e_i$. For each element $\mu_i'$, where $1\leq i \leq \tau$, let $\CB_\G' (\mu_i')=\{ b_j'\, | \, 1\leq j \leq \tau , b_j'\nu = \mu_i'\}$.

 Let us show that always  we can chose  the elements $b_1, \ldots , b_t$ from  the set $\{ b_1',  \ldots , b_\tau'\}$ changing the choice of the elements $a(\g )$ if necessary. Let $$\mu_{i_1}'<\cdots < \mu_{i_\s}'$$ be all the  {\em distinct} elements of the set $\{\mu_1',\ldots , \mu_\tau'\}$.

{\em Definition.} The element $\mu_{i_\alpha}'\in \G $, where $1\leq \alpha \leq \s$,  is called an {\bf avoidable element} if 
$$ \Rel(\mu_{i_\alpha}')\cap \Bigg(\bigcup_{i=1}^s(e_i+\CB_\G )\Bigg)=\emptyset ,$$
and {\bf non-avoidable}, otherwise. The sets of avoidable and non-avoidable elements are denoted by $\G_{av}$ and $\G_{na}$, respectively. Therefore, $$\{\mu_{i_1}',\ldots , \mu_{i_\s}'\}=\G_{av} \coprod \G_{na}.$$ 
{\em Definition.} The set  $\{a(\g ) \, | \, \g \in \G\}$ is called a {\bf non-avoidable} set provided $a(\mu')\not\in \CB_\G' (\mu')$ for all $\mu'\in \G_{na}$. 

It follows from the definition, that there are  plenty of non-avoidable sets. Proposition \ref{A2Jan21} shows that the set $\mB_\G$ is unique provided that $\{a(\g ) \, | \, \g \in \G\}$ is a non-avoidable set.
 
 \begin{proposition}\label{A2Jan21}
We keep the notation as above. Let $\{a(\g ) \, | \, \g \in \G\}$ be a non-avoidable set. Then  the set $\mB_\G$ is uniquely defined. Furthermore,
\begin{enumerate}
\item $\{ b_1, \ldots , b_t\}=\CB_\G'\backslash \{ a(\mu' )\, | \,\mu'\in \G_{av}\}$ and $t=t_\G = |\CB_\G'|-|\G_{av}|$ (see (\ref{BGbas}) and (\ref{BGbas2})). 
\item The set $\{ \mu_{i_1}', \ldots , \mu_{i_\s}'\}$ contains precisely all the distinct elements of the set $\{ \mu_1, \ldots , \mu_t\}$, i.e. $\{ \mu_1, \ldots , \mu_t\}= \{ \mu_{i_1}', \ldots , \mu_{i_\s}'\}$ (after deleting repetitions), and the multiplicity of the element $\mu_i$ is equal to 
$${\rm mult} (\mu_i)=\begin{cases}
|\CB_\G'(\mu_i)|-1& \text{if }\mu_i\in \G_{av},\\
|\CB_\G'(\mu_i)|& \text{if }\mu_i\in \G_{na}.\\
\end{cases}$$
\end{enumerate}
\end{proposition}

{\it Proof}. 1. Since the set $\{a(\g ) \, | \, \g \in \G\}$ is a non-avoidable set, $\{ b_1, \ldots , b_t\}=\CB_\G'\backslash \{ a(\mu' )\, | \,\mu'\in \G_{av}\}$,  and so $t=t_\G = |\CB_\G'|-|\G_{av}|$.

2. Statement 2 follows from statement 1. $\Box $\\

  Summarizing, we have the following proposition which gives generators and defining relations for the algebra $\bAmon$.

\begin{proposition}\label{A24Jun20}
Suppose that $s=|\ind (\G )|\geq 2$ and $ \dec (\G )_{\geq 2}\neq \emptyset$. Then the monomial algebra $\bAmon$ is generated by the minimal set  of generators  $\{ f_{\nu_1}=x^{\nu_1}, \ldots , f_{\nu_s}=x^{\nu_s}\}$  that satisfies the defining relations:
\begin{eqnarray*}
 f^{b_i}&=& f^{a(\mu_i )}\;\; {\rm for   }\;\;i=1, \ldots , t=t_\G ,\\
f^c &=&0\;\; {\rm where  }\;\; c\in \ind \, \N^s(m-1, \G ).
\end{eqnarray*}
None of the relations in the first set can be omitted. 
\end{proposition}

 {\bf Sets of generators and defining relations of the algebra $\bA$ where $A \in \CA (m, \G )$,  $|\ind (\G )|\geq 2$ and $ \dec (\G )_{\geq 2}\neq \emptyset$.} 

\begin{theorem}\label{A6May20}
Given an algebra $A\in \CA (m , \G )$. Suppose that $s=|\ind (\G )|\geq 2$ and $\dec (\G)_{\geq 2}\neq \emptyset$.  Then the algebra $\bA$ is generated by the minimal set of generators    $\{ f_{\nu_1}, \ldots , f_{\nu_s}\}$  that satisfies the defining relations:
\begin{equation}\label{A1st}
f^{b_i}= f^{a(\mu_i)}+\sum_{\g'\in \G (\mu_i )}\th_{\mu_i , \g';b_i}f^{a(\g')}\;\; {\rm for  }\;\; i=1,\ldots , t=t_\G,\;\; {\rm and}
\end{equation}
\begin{equation}\label{A2nd}
f^c =0\;\; {\rm for  }\;\; c\in \ind \, \N^s(m-1, \G ),
\end{equation}
where the elements $\th_{\g_i , \g';*}$ are defined in (\ref{thggb}). None of the relations in (\ref{A1st}) can be omitted.
\end{theorem}

{\it Proof}.  Let $\bA'$ be a commutative algebra that is generated by the indeterminates $f_{\nu_1}, \ldots , f_{\nu_s}$ subject to the defining relations of the theorem. The relations $f^c=0$, where $c\in \ind \, \N^s(m-1, \G ))$, imply that the algebra $\bA'$ is a finite dimensional algebra. The $K$-homomorphism
$$ \v : \bA'\ra \bA , \;\; f_{\nu_i}\mapsto f_{\nu_i}$$ is an epimorphism. To finish the proof it suffices to show that $\dim_K(\bA')\leq \dim_K(A)$. The algebra $\bA'$ admits a descending  filtration $\{ \bA'_{\geq i}:=\sum_{\{a\in \N^s\, | \, a\nu\geq i\}} Kf^a\}_{i\in \N}$. Its image under $\v$, $\{ \v (\bA'_{\geq i})\}_{i\in \N}$, is the induced radical filtration on the algebra $\bA$. Hence, we have the epimorphism of the corresponding graded algebras:
$$ {\rm gr}\, \v : {\rm gr}\, \bA'\ra {\rm gr}\, \bA =K\oplus \bigoplus_{\g \in \G} Kx^\g =\bAmon .$$
  The algebra ${\rm gr}\, \bA'$ is an epimorphic image of the $\N$-graded algebra $\L$ which is generated by the indeterminates $f_{\nu_1}, \ldots , f_{\nu_s}$ subject to the defining relations $f^{b_i}= f^{a(\mu_i)}$ for $i=1, \ldots ,t$ and $f^c=0$ for all $c\in \ind \, \N^s(m-1,\G )$.  By Proposition \ref{A24Jun20}, $\L \simeq \bAmon$. Hence, $\dim_K(\bA')\leq \dim_K(\bA)$, as required. $\Box$\\





{\bf Generators and defining relations for the algebra $A\in \CA (m, \G )$.} For each algebra $A\in \CA (m, \G )$, we give explicit sets of generators and defining relations using the sets of generators and defining relations of the algebra $\bA$ that we obtained above. 

{\it Case:} $\G = \emptyset$. The algebra $\L = \L (M):=K\oplus x^mK[x]$ is the only element of the set $\CA (m, \emptyset )$. Lemma \ref{d7May20} follows at once from the equality $\L(m) = \bigoplus_{i=0}^{m-1}K[x^m]x^{m+i}$, a direct sum of free rank 1 $K[x^m]$-modules.

\begin{lemma}\label{d7May20}
The commutatiove algebra $\L (m)$ is generated by the elements $\{ x^{m+i}\, | \, i=0,1,\ldots , m-1\}$ subject to the defining relations: For all $i,j$ such that $0\leq i,j\leq m-1$, 
$$x^{m+i}x^{m+j}=\begin{cases}
x^mx^{m+i+j}& \text{if }i+j<m,\\
(x^m)^2x^{i+j}& \text{if }i+j\geq m.\\
\end{cases}$$
The set  $\{ x^{m+i}\, | \, i=0,1,\ldots , m-1\}$ is a minimal set of generators and the defining relations above is a minimal set of defining relations. 
\end{lemma}

We have the direct sum of vector spaces over $K$, $K[x]=K[x]_{\leq m-1}\oplus (x^m)$ where $K[x]_{\leq m-1}=\bigoplus_{i=0}^{m-1}Kx^i$. We identify the factor algebra $$F=K[x]/(x^m)=\bigoplus_{i=0}^{m-1}Kx^i$$ with the $K$-subspace 
 $K[x]_{\leq m-1}$ via the $K$-linear isomorphism $x^i\mapsto x^i$ for $i=0,1,\ldots , m-1$. Then every polynomial $p\in K[x]$ is a {\em unique} sum 
 $$p=\bp+[p]\;\; {\rm where}\;\;\bp \equiv p\mod (x^m)\;\; {\rm  and }\;\; [p]=p-\bp\in (x^m).$$ In particular, for elements $f,g\in F$, there are two products: $fg$ is in $K[x]$ and $f\cdot g$ is in $F$. They are related by the equality in $K[x]$:
 $$f\cdot g = fg-[fg].$$
 We will use this type of equalities when we obtain defining relations for the algebra $A$ from the defining relations of the algebra $\bA$.  Recall that 
 $$\L (m) =\bigoplus_{i=0}^{m-1} K[x^m]x^{m+i}\subseteq K[x].$$
 So, for each polynomial $p\in K[x]$, its  projection $[p]$ onto $(x^m)$ is a unique sum
\begin{equation}\label{pbpb}
[p]=\sum_{i=0}^{m-1} p_i(x^m)x^{m+i}\;\; {\rm where}\;\; p_i(x^m)\in K[x^m].
\end{equation}
 When we write $[p]$ we mean that the element $[p]\in (x^m)$ is written as the unique sum above.

\begin{theorem}\label{Aa5May20}
We keep the notation of Corollary \ref{a5May20}. 
 Given an algbera $A\in \CA (m, \G )$ with $\G\neq \emptyset$.  Then the commutative algebra $A$ is generated by the elements $\{ f_\g, x^{m+i}\, | \,\g \in \G , i=0,1,\ldots , m-1\}$ subject to the defining relations:
\begin{eqnarray*}
 f_\g f_{\g'}-[f_\g f_{\g'}]&=& f_{\g+\g'}+\sum_{\rho \in \G (\g +\g')}\mu_{\g, \g';\rho}f_\rho\;\; {\rm for\; all}\;\; \g , \g'\in \G\;\; {\rm such \; that}\;\; \g+\g'<m-1,\\
 f_\g f_{\g'}-[f_\g f_{\g'}]&=& 0\;\; {\rm for\; all}\;\; \g , \g'\in \G\;\; {\rm such \; that}\;\; \g+\g'\geq m-1,\\
 x^{m+i}x^{m+j}&=& [x^{2m+i+j}]\;\; {\rm for\; all}\;\; i,j=0,1,\ldots , m-1, \\
 x^{m+i}f_\g &=& [ x^{m+i}f_\g]\;\; {\rm for\; all}\;\; i,j=0,1,\ldots , m-1 \;\; {\rm and }\;\; \g \in \G .
\end{eqnarray*}
\end{theorem}

{\it Proof}. Let $A'$ be the algebra which is defined by the generators and defining relations  in the theorem. We have to show that it is isomorphic to the algebra $A$.  By Corollary  \ref{a5May20}, the   algebra $A$ is generated by the elements 
$\{ f_\g, x^{m+i}\, | \,\g \in \G , i=0,1,\ldots , m-1\}$ that satisfy the relations of theorem. Hence, there is a a natural $K$-algebra epimorphism: 
$$\v : A'\ra A, \;\; f_\g \mapsto f_\g , \;\; x^{m+i}\mapsto x^{m+i}.$$
By Lemma \ref{d7May20}, the algebra $\L (m)$ is a subalgebra of $A'$ that is mapped isomorphically via $\v$ onto its copy $\L (m)$ in the algebra $A$. The ideal $\ga =\sum_{i=0}^{m-1}K[x^m]x^{m+i}$ of the subalgebtra $\L (m)$ of $A'$ is an ideal of the algebra $A'$ (see the last two types of the defining relations of the algebra $A'$). It is mapped isomorphically by $\v$ onto the ideal $(x^m)=x^mK[x]$ of the subalgebra $\L (m)$ of $K[x]$. By Corollary \ref{a5May20}, $\bA':=A'/\ga \simeq \bA$ (as the algebra $\bA'$ has the same generators  and defining relations as the algebra $\bA$). Now, there is a commutative diagram of algebra homomorphisms:

\begin{displaymath}
    \xymatrix{
 0\ar[r]&  \ga\ar[r]\ar[d]&       A' \ar[r] \ar[d]_\v & \bA'  \ar[r]\ar[d] &0\\
  0\ar[r]& (x^m)\ar[r]&     A \ar[r]       & \bA \ar[r] &0}
\end{displaymath}

where the left  vertical map is a restriction of $\v$ and the right   vertical map is induced by $\v$, they are bijections. Hence, the map $\v$ is so, i.e. $A'\simeq A$ via $\v$. $\Box $\\

{\it Definition.} For the algebra $A\in \CA (m, \G )$, the $K$-basis $\{ 1, f_\g , x^i\, | \, \g \in \G, i\geq m$ is called {\em the canonical basis} of $A$. By definition, {\em a canonical basis} of $A$ is obtained from the canonical basis my multiplying each its element by an element of $K^\times$. \\

{\bf The structure constants of the algebra $A\in \CA (m , \G )$.} The structure constants of the algebra $A\in \CA (m , \G )$ with respect to  the canonical basis  is given by the rule: For all elements $\g, \g'\in \G$ and $i, j\geq m$,

\begin{eqnarray*}
 f_\g f_{\g'}&=& \begin{cases}
f_{\g +\g'}+\sum_{\rho \in \G (\g +\g')}\mu_{\g , \g'; \rho}f_\rho+[f_\g f_{\g'}]& \text{if }\g+\g'<m-1,\\
[f_\g f_{\g'}]& \text{if }\g+\g'\geq m-1,\\
\end{cases}\\
x^if_\g &=&[x^if_\g ]_m,\\
x^ix^j&=& x^{i+j},
\end{eqnarray*}

where the element $[\cdots ]_m$ has to be  written as a sum $\sum_{j\geq m}\l_jx^j$ with $\l_j\in K$. \\

Lemma \ref{Aa7May20} gives generators and defining relations for algebras $A\in \CA (m, \G )$ such that  $|\ind (\G )| =1$. 

\begin{lemma}\label{Aa7May20}
Given an algebra $A\in \CA (m, \G )$. Suppose that $\ind (\G ) =\{ \nu_1\}$. Then the algebra $A$ is generated by the elements $\{ f_{\nu_1} , x^{m+i}\, | \, i=0,1,\ldots , m-1\}$ subject to the defining relations: 
\begin{eqnarray*}
 f_{\nu_1}^{[\frac{m-1}{\nu_1}]+1}&=&\Bigg[f_{\nu_1}^{[\frac{m-1}{\nu_1}]+1}\Bigg], \\
x^{m+i}x^{m+j} &=&[x^{2m+i+j}]\;\;{\rm for}\;\; i,j=0,1,\ldots , m-1,\\
x^{m+i}f_{\nu_1}&=& [x^{m+i}f_{\nu_1}]\;\;{\rm for}\;\; i=0,1,\ldots , m-1.
\end{eqnarray*}
\end{lemma}

{\it Proof}. The lemma follows from Lemma \ref{a7May20}. $\Box$ \\

Theorem  \ref{X6May20} gives generators and defining relations for algebras $A\in \CA (m, \G )$ such that  $|\ind (\G )| \geq 2$. In general, some of the relations that are given in the first equality of Theorem \ref{X6May20} are redundant. An irredundant set is given in Theorem \ref{XA6May20}.

\begin{theorem}\label{X6May20}
We keep the notation of Theorem \ref{6May20}. 
 Given an algebra $A\in \CA (m , \G )$. Suppose that  $s=|\ind (\G )|\geq 2$. Then the algebra $A$ is generated by the elements   $\{ f_{\nu_1}, \ldots , f_{\nu_s}, x^{m+i}\, | \, i=0,1,\ldots , m-1\}$ subject to the defining relations (as a commutative algebra):
\begin{eqnarray*}
 f^b-[f^b]&=& f^{a(\g )}-[f^{a(\g )}]+\sum_{\g'\in \G (\g )}\th_{\g , \g';b}(f^{a(\g')}-[f^{a(\g')}])\; {\rm where  }\; \g \in \dec (\G )_{\geq 2}\; {\rm and}\; b\in \Rel(\g )\backslash \{ a(\g )\},\\
f^c &=&[f^c]\;\; {\rm where  }\;\; c\in \ind \, \N^s(m-1, \G ),\\
x^{m+i}f_{\nu_j}&=& [x^{m+i}f_{\nu_j}]\;\;{\rm for}\;\; i=0,1,\ldots , m-1\;\; {\rm and}\;\; j=1,\ldots , s, \\
x^{m+i}x^{m+j} &=&[x^{2m+i+j}]\;\;{\rm for}\;\; i,j=0,1,\ldots , m-1.
\end{eqnarray*}
\end{theorem}

{\it Proof}. Repeat the proof of Theorem \ref{Aa5May20} word for word but the algebra $A'$ there is replaced by the algebra $A'$ from the theorem and instead of using Corollary \ref{a5May20} to prove the isomorphism $A'\simeq A$ we use Theorem \ref{6May20} instead. $\Box$\\

Corollary \ref{Ab7May20} is a particular (degenerated) case of Theorem \ref{X6May20}.
 For each algebra $A\in \CA (m, \G )$ such that  $|\ind (\G )| \geq 2$ and $\dec (\G )_{\geq 2}=\emptyset$ it  gives generators and defining relations.
 
\begin{corollary}\label{Ab7May20}
 Given an algebra $A\in \CA (m , \G )$. Suppose that  $s=|\ind (\G )|\geq 2$ and $\dec (\G )_{\geq 2}=\emptyset$. Then the algebra $A$ is generated by the elements   $\{ f_{\nu_1}, \ldots , f_{\nu_s}, x^{m+i}\, | \, i=0,1,\ldots , m-1\}$ subject to the defining relations (as a commutative algebra):
\begin{eqnarray*}
f^c &=&[f^c]\;\; {\rm where  }\;\; c\in \ind \, \N^s(m-1, \G ),\\
x^{m+i}f_{\nu_j}&=& [x^{m+i}f_{\nu_j}]\;\;{\rm for}\;\; i=0,1,\ldots , m-1\;\; {\rm and}\;\; j=1,\ldots , s, \\
x^{m+i}x^{m+j} &=&[x^{2m+i+j}]\;\;{\rm for}\;\; i,j=0,1,\ldots , m-1.
\end{eqnarray*}
\end{corollary}

Theorem  \ref{XA6May20} gives generators and defining relations for algebras $A\in \CA (m, \G )$ such that  $|\ind (\G )| \geq 2$ and $\dec (\G)_{\geq 2}\neq \emptyset$, 
and none of the relations given in the first equality of Theorem \ref{XA6May20} is redundant.

\begin{theorem}\label{XA6May20}
We keep the notation of Theorem \ref{A6May20}. Given an algebra $A\in \CA (m , \G )$. Suppose that $s=|\ind (\G )\geq 2$ and $\dec (\G)_{\geq 2}\neq \emptyset$.  Then the algebra $A$ is generated by the  set   $\{ f_{\nu_1}, \ldots , f_{\nu_s}, x^{m+i}\, | \, i=0,1,\ldots , m-1\}$ subject to the defining relations:
\begin{eqnarray*}
f^{b_i}-[f^{b_i}] &=& f^{a(\mu_i)}-[f^{a(\mu_i)}]+\sum_{\g'\in \G (\mu_i )}\th_{\mu_i , \g';b_i}(f^{a(\g')}-[f^{a(\g')}])\;\; {\rm for  }\;\; i=1,\ldots , t=t_\G,\\f^c &=&[f^c]\;\; {\rm where  }\;\; c\in \ind \, \N^s(m-1, \G ),\\
x^{m+i}f_{\nu_j}&=& [x^{m+i}f_{\nu_j}]\;\;{\rm for}\;\; i=0,1,\ldots , m-1\;\; {\rm and}\;\; j=1,\ldots , s, \\
x^{m+i}x^{m+j} &=&[x^{2m+i+j}]\;\;{\rm for}\;\; i,j=0,1,\ldots , m-1.
\end{eqnarray*}
None of the relations in the first type of relations is redundant.
\end{theorem}

{\it Proof}.  Repeat the proof of Theorem \ref{Aa5May20} word for word but the algebra $A'$ there is replaced by the algebra $A'$ from the theorem and instead of using Corollary \ref{a5May20} to prove the isomorphism $A'\simeq A$ we use Theorem \ref{A6May20} instead. $\Box$


\section{Isomorphism problems and an explicit description of the automorphism groups of algebras   $A\in \CA (m)$}\label{EDAUTA}

{\em In this section, the field $K$ is an algebraically closed field (unless it is not stated otherwise)}.\\

The aim of the section is to give explicit generators for the automorphism groups of algebras in $\CA (m)$ (Teorem \ref{28Apr20}.(1)), to give criteria for the automorphism group  to be infinite/finite/trivial (Theorem \ref{28Apr20}.(2--4), and Corollary  \ref{a2May20}). Theorem \ref{2May20} and Corollary \ref{c3May20} are explicit descriptions of the  orders of all finite automorphism groups of algebras  in $\CA (m)$.

Let $\mO (m)$ be the set of orders of all {\em finite} automorphism groups of algebras in $\CA (m)$. By Corollary \ref{a30Apr20}.(1),
$$ \mO (m) \subseteq \begin{cases}
\{ 1, \ldots , m-3\}& \text{if } 4\leq m \text{ is even},\\
\{ 1, \ldots , m-4\}& \text{if } 5\leq m \text{ is odd}.\\
\end{cases}
$$
Theorem \ref{2May20} and Corollary \ref{c3May20}
 are explicit descriptions of the set $\mO (m)$. Lemma \ref{a23Apr20} and Lemma \ref{d3May20} are examples of algebras in $\CA (m)$ with sophisticated orders of automorphism groups. 
 
 Notice that $\Aut_K(K[x])=\{ \s_{\l\mu}\, | \, \l \in K^\times, \mu \in K\}$ where $\s_{\l \mu} (x) = \l  x+\mu$. The group $$\Aut_K(K[x])=\Sh (K)\rtimes \mT (K)$$ is a semi-direct product of its normal subgroup, the {\em shift group},  $\Sh (K):=\{ s_\mu :=\s_{1,\mu}\, | \, \mu \in K\} \simeq (K, +)$ ($s_\mu \mapsto \mu$) and the {\em algebraic torus} $\mT = \mT (K):=\{ t_\l :=\s_{\l , 0}\, | \, \l\in K^\times \} \simeq (K^\times , \cdot )$ ($t_\l\mapsto \l$). The automorphism group $\Aut_K(K[x])$ acts in the obvious way on the set of ideals of the polynomial algebra $K[x]$. Then the algebraic torus $\mT$ is the stabilizer  of each of the ideals $(x^i)$, $i\geq 1$, i.e. 
$\mT = \{ \s \in \Aut_K(K[x])\, | \, \s ((x^i))= (x^i)\}$. 

For a natural number $n\geq 2$, an element $\l \in K^\times$ is called a {\em primitive $n$'th root of unity} if $\l^n=1$ and the elements $1, \l , \ldots , \l^{n-1}$ are distinct. If ${\rm char }(K)=0$ then for each $n\geq 2$ there are primitive $n$'th roots of unity. If  ${\rm char }(K)=p>0$  then for each  natural number $n\geq 2$ there are primitive $n$'th roots of unity  iff  $p\nmid n$ (the field $K$ is algebraically closed). If  ${\rm char }(K)=0$   (resp., ${\rm char }(K)=p>0$) then for each natural number $n\geq 2$ (resp., such that $p\nmid n$) fix a primitive $n$'th root of unit, say $\l_n$. Then the identity group and  {\em the cyclic groups of order} $n\geq 2$ if ${\rm char} (K)=0$ (resp., ${\rm char} (K)=p>0$ and $p\nmid n$), 
\begin{equation}\label{Cndef}
C_n:=\langle t_{\l_n}\rangle= \{ t_{\l_n}^i\, | \, 0\leq i \leq n-1\},
\end{equation}
{\em are precisely all the finite subgroups of the algebraic torus} $\mT$.\\

{\bf The set $\Iso_K(A,A')$ where $A\in \CA (m)$ and $A'\in \CA (m')$.} For an algebra $A\in \CA (m)$,  let $\{ 1, f_\g=x^\g u_\g\, | \, \g \in \G_A\}$ be the  canonical basis of the algebra $\bA$ where for each $\g \in \G$,  $f_\g = x^\g+\sum_{\d\in C\G (\g )}\l_{\g \d}x^{\d}=x^\g u_\g$ (see (\ref{fgg})) where  
\begin{equation}\label{ugdef}
u_\g :=1+ \sum_{\d\in C\G (\g )}\l_{\g \d}x^{\d-\g}\in \bA^\times 
\end{equation}
and $\bA^\times $ is the group of units of the algebra $\bA$. 

{\em Definition.} The elements $\{ u_\g \, | \, \g \in \G_A\}$ of $\bA^\times $ are called the {\bf canonical units} of the algebra $\bA$.

Theorem \ref{3May20}.(2) is an isomorphism criterion for algebras in $\bigcup_{m\geq 2}\CA (m)$.   Theorem \ref{3May20}.(3) describes the set $\Iso_K(A,A')$ of all $K$-algebra isomorphisms from an algebra $A$ to $A'$. 

\begin{theorem}\label{3May20}
($K$ is an arbitrary field) Given algebras $A\in \CA (m)$ and $A'\in \CA (m')$. Let $\{ 1, f_\g=x^\g u_\g\, | \, \g \in \G_A\}$ and $\{ 1, f_{\g'}'=x^{\g'} u_{\g'}\, | \, \g' \in \G_{A'}\}$ be the canonical bases of the algebras $A$ and $A'$, respectively. Then 
\begin{enumerate}
\item $\Iso_K(A,A')=\{ t_\l\in \mT \, | \, t_\l(A)=A'\}$. 
\item $A\simeq A'$ iff $m=m'$, $\G_A=\G_{A'}$ and there exists an automorphism $t_\l \in \mT$ such that $t_\l (f_\g)=\l^\g f_{\g}'$ for all $\g \in \G_A$ iff $m=m'$, $\G_A=\G_{A'}$ and there exists an automorphism $t_\l \in \mT$ such that $t_\l (u_\g)= u_{\g}'$ for all $\g \in \G_A$.
\item Suppose that $A\simeq A'$. Then 
\begin{enumerate}
\item $\Iso_K(A,A')=\{ t_\l\in \mT \, | \, t_\l(f_\g)=\l^\g f_{\g}'$ for all $\g \in \G_A=\G_{A'}\}=\{ t_\l\in \mT \, | \, t_\l(u_\g)=u_{\g}'$ for all $\g \in \G_A=\G_{A'}\}$. 
\item The $\mT$-orbit $\mT A$ of the algebra $A$
 (i.e. the set of all algebras isomorphic to $A$, by statement 1) is equal to the set $\mT / \Aut_K(A)=\{ t_\l \Aut_K(A)\, | \, t_\l \in \mT\}$ of all subalgebras of $K[x]$ that are isomorphic to the algebra $A$ where $\Aut_K(A)=\{ t_\l \in \mT \, | \, t_\l (u_\g ) = u_\g$ for all $\g \in \G_A\}$. 
\end{enumerate}

\end{enumerate}
\end{theorem}

{\it Proof}. 1. We can assume that $A, A' \subseteq K[x]$ (since the polynomial algebra $K[x]$ is their commom normalization). Suppose that $\s : A\ra A'$ is a $K$-isomorphism. Then it can be uniquely extended to an automorphism $\s : K[x]\ra K[x]$, and necessarily $\s ((x^m))=((x^{m'}))$ (since $(x^m) $ and $(x^{m'})$ are the conductors of the algebras $A$ and $A'$, respectively). Then $m=m'$ and $\s (x) = \l  x$ for some $\l \in K^\times$, i.e. $t_\l \in \mT$, and statement 1 follows.

2. It suffices to prove that the first `iff' holds as the second one follows from the first.

$(\Rightarrow )$ Suppose that $A\simeq A'$. By statement 1, there is an automorphism $t_\l \in \mT$ such that $t_\l (A) = A'$. Then $m=m'$ (see the proof of statement 1) and the automorphism $t_\l$ of $K[x]$ respects the $(x)$-adic filtration of the algebra $K[x]$. As a result, it respects the induced $(x)$-adic filtrations on the algebras $\bA =A/x^mK[x]$ and $\bA' =A'/x^mK[x]$. Hence, the automorphism $t_\l$ induces an automorphism of the associated graded algebras, ${\rm gr} (\bA)\simeq {\rm gr} (\bA')$. Therefore, $\G_A = \G_{A'}$ and the image of {\em the } canonical basis of the algebra $A$ under $t_\l$  is {\em a } canonical basis of the algebra $\bA'$, and so $t_\l (f_\g ) = \l^\g f_\g'$ for all elements $\g \in \G_A=\G_{A'}$. 

$(\Leftarrow )$ Clearly, $t_\l (A) =A'$, and so $t_\l \in \Iso_K(A, A')$, by statement 1. 

3(a). The statement (a) follows from statement 2.

(b). The statement (b) follows from the statement (a) and the fact that $\Aut_K(A) = \{ t_\l \in \mT \, | \, t_\l (u_\g ) = u_\g$ for all $\g \in \G_A\}$, by the statement (a) where $A=A'$. 
$\Box $\\

 {\bf Description of orbits of the algebraic torus $\mT$ action on  $K^{\times s}$.} For each natural number $n\geq 1$, let $\mU_n=\mU_n(K):=\{\l \in K\, |\, \l^n=1\}$ be the group of $n$'th roots of unity. In characteristic zero the set $\mU_n$ is a cyclic group of order $n$. In prime characteristic $p>0$, $\mU_n=\mU_{n'}$ is a cyclic group of order $n'$ where $n=p^\nu n'$ for unique natural numbers $\nu$ and $n'$ such that $(p,n')=1$.  A cyclic generator of the group $\mU_n$ is called a primitive $n$'th root of unity if char$(K)=0$ and 
 a primitive $n'$'th root of unity if char$(K)=p>0$.
 
  Let us consider an action of the algebraic torus $\mT$ on the set $K^{\times s}$, $s\geq 1$: For all $\l\in\mT$ and $ (\l_1, \ldots , \l_s)\in K^{\times s}$, 
\begin{equation}\label{TKs}
\l\cdot (\l_1, \ldots , \l_s)=(\l^{n_1}\l_1, \ldots , \l^{n_s}\l_s),
\end{equation} 
 where the integer-valued vector $(n_1,\ldots , n_s)\in (\Z\backslash \{ 0 \})^s$ is called the {\em  weight vector} of the action. For simplicity, we consider the case when all the weights $n_i$ are not equal to zero. The general case is easily reduced to this one. 
 
   Let $\xi_1$ be a cyclic generator of the group $\mU_{n_1}$, i.e. $\xi_1$ is a primitive $n_1$'th (resp., $n_1'$'th)  root of unity if char$(K)=0$ (resp., char$(K)=p>0$, $n_1=p^{\nu_1}n_1'$ and $p\nmid n_1$).  Let $H(n_1, \ldots , n_s)$ be a (cyclic) subgroup of the group $\mU_{n_1}^{s-1}$ which is  generated by the 
 element $(\xi_1^{n_2}, \ldots ,\xi_1^{n_s})$. 
The order of the group $H(n_1, \ldots , n_s)$ is equal to $n_1\gcd (n_1, \ldots , n_s)^{-1}$ (resp., $n_1'\gcd (n_1',\ldots , n_s')^{-1}$) if 
  char$(K)=0$ (resp., char$(K)=p>0$) where $n_i=p^{\nu_i}n_i'$ and $p\nmid n_i'$.
 
Proposition \ref{A19Jul20} is a criterion for two elements of the set $K^{\times s}$ to belong to the same $\mT$-orbit.
 
 \begin{proposition}\label{A19Jul20}
Given an action of the algebraic torus $\mT$ on the set $K^{\times s}$ with weight vector $(n_1,\ldots , n_s)\in (\Z\backslash \{ 0 \})^s$.  Then elements $\l , \l'\in K^{\times s}$ belong to the same $\mT$-orbit iff 
$ \mu (\l' )\mu (\l )^{-1}\in
H(n_1, \ldots , n_s)$ where for $\l =(\l_1, \ldots , \l_n)$, $\mu (\l )=(\mu (\l )_2,\ldots , \mu(\l)_n)$ and $\mu (\l )_i=\l_1^{-\frac{n_i}{n_1}}\l_i$.
\end{proposition}

{\it Proof}. Notice that $\mT\l = \mT (1,\mu (\l))$ and $\mT\l' = \mT (1,\mu (\l'))$. Hence, $\mT\l' = \mT\l$ iff $\mu (\l')=\mu (\l) (\xi_1^{n_2}, \ldots ,\xi_1^{n_s})^i$ for some natural number $i\geq 0$, and the result follows. $\Box $

Let $A\in \CA (m,\G)$ and $\ind (\G )=\{\nu_1,\ldots , \nu_s\}$. Recall that $f_{\nu_i}=x^{\nu_i}u_{\nu_i}$ where $u_{\nu_i}=1+\sum_{j\in C\G (\nu_i)}\l_{\nu_i,j}x^{j-\nu_i}$ (see (\ref{ugdef})), and for all $\l\in K^\times$, 
$$t_\l (u_{\nu_i})=1+\sum_{j\in C\G (\nu_i)}\l^{j-\nu_i}\l_{\nu_i,j}x^{j-\nu_i}.$$
For each $i=1, \ldots, s$, let $c(u_{\nu_i})=(\ldots, \l_{\nu_i, j}, \ldots)$, $j\in C\G (\nu_i)$,  be the vector of nonzero coefficients 
 (excluding 1)  of the polynomial $u_{\nu_i}$ and let $n(u_{\nu_i})= (\ldots, j-\nu_i, \ldots)$, $j\in C\G (\nu_i)$,  be the corresponding weight vector. Let $c(A)=(c(u_{\nu_1}), \ldots , c(u_{\nu_s}))$ and $n(a)=(n(u_{\nu_1}), \ldots , n(u_{\nu_s}))$.

\begin{proposition}\label{a20Jul20}
($K$ is an arbitrary field) Let $A\in \CA (m,\G)$, $A'\in \CA (m',\G')$ and $\ind (\G )=\{\nu_1,\ldots , \nu_s\}$. Then $A\simeq A'$ iff $\G = \G'$, $m=m'$, $n(u_{\nu_i})=n(u_{\nu_i}')$ for all $i=1,\ldots , s$ and $\mT c(A)=\mT c(A')$.
\end{proposition}

{\it Proof}. The result follows at once from Theorem \ref{3May20}.(2). $\Box $\\

{\em Remark.} Using  Proposition \ref{A19Jul20} and Proposition \ref{a20Jul20}, we can establish in finitely mant steps whether algebras $A, A'\in \CA$ are isomorphic or not.

{\em Definition.} Let $f(x) = x^d+a_{d-1}x^{d-1}+\cdots+a_1x+a_0\in K [x]$ be a monic polynomial of degree $d\geq 1$ where $a_i\in K$ are the coefficients of the polynomial $f(x)$. Then the natural number
$$ \gcd (f(x)):=\gcd \{ i\geq 1\, | \, a_i\neq 0\}$$
is called the {\em exponent} of $f(x)$.

Clearly, the exponent of $f(x)$ is the largest natural number $m\geq 1$ such that $f(x)=g(x^m)$ for some polynomial $g(x)\in K[x]$. \\

{\it Definition.} Suppose  that ${\rm char} (K)=p>0$. Every nonzero natural number $n$ is a unique product of natural numbers,  $n=p^sn_p$ where $s\geq 0$ and $p\nmid n_p$. The number $n_p$ is called the $p$-{\em co-prime divisor} of $n$. In particular, for a nonscalar polynomial $f\in K[x]$, we have the  unique product
\begin{equation}\label{gcdpf}
 \gcd (f) = p^s  \gcd_p (f)\;\; {\rm  where }\;\; s\geq 0, \;\; \gcd_p (f)\in \N\;\; {\rm  and}\;\;  p\nmid
\gcd_p (f).
\end{equation}
Let $\{ 1, f_\g=x^\g u_\g\, | \, \g \in \G =\G_A\}$ be the canonical basis of the algebra $\bA$, Let 
\begin{eqnarray*}
\gcd (\bA ) &:=& \gcd \{ \gcd (u_\g )\, | \, \g \in \G\}, \\
\gcd_p (\bA ) &:=& \gcd \{ \gcd_p (u_\g )\, | \, \g \in \G\}, 
\end{eqnarray*}
where $p$ is a prime number. Clealry, $\gcd_p(\bA ) = \gcd (\bA )_p$, the $p$-co-prime divisor of $\gcd (A)$. \\

{\bf The group $\Aut_K(A)$ where $A\in \CA (m)$.}  Theorem \ref{28Apr20} is an explicit description of the automorphism groups of  algebras in $\CA (m)$. 
\begin{theorem}\label{28Apr20}
Let $A\in \CA (m, \G )$, $\{ 1, f_\g =x^\g u_\g\, | \, \g \in \G\}$ be the canonical basis of the algebra $\bA$ and $\ind (\G ) = \{ \nu_1, \ldots , \nu_s\}$. Then 
\begin{enumerate}
\item $\Aut_K(A)=\{ t_\l \in \mT \, | \, t_\l (f_\g ) = \l^\g f_\g$ for all $\g \in \G\}=\{ t_\l \in \mT \, | \, t_\l (u_\g ) = u_\g$ for all $\g \in \G\}$. 
\item $\Aut_K(A) =\mT$ iff $f_\g =x^\g$ for all $\g \in \G$.
\item Suppose that $f_\g \neq x^\g$ for some  $\g \in \G$ (i.e. $\Aut_K(A) \neq \mT$, by statement 2). Then $\Aut_K(A)=C_n=\langle t_{\l_n}\rangle$ is a cyclic group of order $n$ where 
$$ n= \begin{cases}
\gcd (A)& \text{if char}(K)=0,\\
\gcd (A)_p& \text{if char}(K)=p>0,\\
\end{cases}
$$
and $\l_n$ is a primitive $n$'th root of unity.
\begin{enumerate}
\item If ${\rm char}(K)=0$ then $\gcd (A) = \gcd \{ \gcd (u_{\nu_i})\, | \, i=1, \ldots , s\}$.
\item If ${\rm char}(K)=p>0$ then $\gcd (A)_p = \gcd \{ \gcd_p (u_{\nu_i})\, | \, i=1, \ldots , s\}$.
\end{enumerate}
\item $\Aut_K(A)=\{ e\}$ iff $\gcd (A) =1$ if ${\rm char}(K)=0$, or $\gcd_p (A) =1$ if ${\rm char}(K)=p>0$.
\end{enumerate}
\end{theorem}

{\it Proof}. 1. Statement 1 follows from Theorem \ref{3May20}.(3).

2. Statement 2 follows from statement 1. 

3. The first statement of statement 3 follows from statement 1. 

(a,b). Let 
$$ n= \begin{cases}
\gcd (A)& \text{if char}(K)=0,\\
\gcd (A)_p& \text{if char}(K)=p>0,\\
\end{cases}\;\; {\rm and}\;\; 
m= \begin{cases}
\gcd \{ \gcd (u_{\nu_i})\, | \, i=1, \ldots , s\}& \text{if char}(K)=0,\\
\gcd \{ \gcd_p (u_{\nu_i})\, | \, i=1, \ldots , s\}& \text{if char}(K)=p>0.\\
\end{cases}
$$
Then $n|m$. To finish the proof it is enough to show that $ m|n$. Let $t:=
t_{\l_m}$. Then $t\in \mT$ and 
 $t(f_{\nu_i})\in K^\times f_{\nu_i}$ for all $i=1, \ldots , s$. The elements $f_{\nu_1}, \ldots , f_{\nu_s}$ are generators of the algebra $A$. Hence, $t(A)=A$, and so $t\in \Aut_K(A)$. The order of the automorphism $t$, which is $m$, divides the order of the group $\Aut_K(A)$, which is $n$, as required. 

4. Statement 4 follows fro statement 3. 
$\Box $\\

 Corollary \ref{a2May20} is a criterion for an algebra $A \in \CA (m)$ to have an infinite automorphism group.

\begin{corollary}\label{a2May20}
Given an algebra $A \in \CA (m, \G)$. Then the algebra $A$ is a monomial algebra iff $\Aut_K(A)=\mT$ iff $\{  1, x^\g \, | \, \g \in \G\}$ is the canonical basis of the algebra $\bA$. 
\end{corollary}

{\it Proof}. The corollary follows from Theorem \ref{28Apr20}. $\Box $\\

Corollary \ref{a30Apr20} gives the upper bound for the orders of finite automorphism groups of algebras in $\CA (m)$.

\begin{corollary}\label{a30Apr20}
 
\begin{enumerate}
\item Suppose that $A\in \CA (m)$ and $|\Aut_K(A)|<\infty$ (i.e. the algebra $A\in \CA (m)$ is not the  monomial algebra). Then 
$$ |\Aut_K(A)|\leq \begin{cases}
m-3& \text{if } 4\leq m \text{ is even},\\
 m-4& \text{if } 5\leq m \text{ is odd}.\\
\end{cases}
$$
\begin{enumerate}
\item If char$(K)=0$ then the upper bounds above are the exact upper bounds, see statements 2(a,b) below.
\item If char$(K)=p>0$ then 
$$ |\Aut_K(A)|\leq \begin{cases}
\max \{ i_p \, | \, 1\leq i \leq m-3 \}& \text{if } 4\leq m \text{ is even},\\
 \max \{ i_p \, | \, 1\leq i \leq m-4 \}& \text{if } 5\leq m \text{ is odd},\\
\end{cases}
$$
where $i_p$ is the $p$-co-prime divisor of the natural number $i$. Let $u(p,m):=\max\{i_p\, | \, 1\leq i\leq m-3$, $i+j\leq m-1$ for some natural number $j=j(i)\geq 2$ such that $j\nmid m-1\}$. Then $u(p,m)=\max\{ |\Aut_K(A)|\; | \; A\in \CA (m), |\Aut_K(A)|<\infty\}$.
\end{enumerate}
\item
\begin{enumerate}
\item If $m\geq 4$ is even and $A=\sum_{i\geq 0} Kg^i+x^mK[x]$ where $g=x^2(1+x^{m-3})$ then $\G_A= \{ 2i\, | \, 1\leq i\leq \frac{m-1}{2}\}$, $\{ 1, g, x^{2i}\, | \, 2\leq i \leq \frac{m-1}{2}\}$ is the canonical basis of the algebra $\bA$ and $$\Aut_K(A)=\begin{cases}
\langle t_{\l_{m-3}}\rangle & \text{if char}(K)=0,\\
\langle t_{\l_{(m-3)_p}}\rangle& \text{if char}(K)=p>0,\\
\end{cases} 
|\Aut_K(A)|=\begin{cases}
m-3 & \text{if char}(K)=0,\\
(m-3)_p& \text{if char}(K)=p>0,\\
\end{cases}
$$ where $\l_{m-3}$ is a primitive $(m-3)$'rd root of unity, $(m-3)_p$ is the $p$-co-prime divisor of $m-3$ if ${\rm char} (K)=p>0$. 
\item  If $m\geq 5$ is odd, $3\nmid m-1$  and $A=\sum_{i\geq 0} Kg^i+x^mK[x]$ where $g=x^3(1+x^{m-4})$ then $\G_A= \{ 3i\, | \, 1\leq i\leq \frac{m-1}{3}\}$, $\{ 1, g, x^{3i}\, | \, 2\leq i \leq \frac{m-1}{3}\}$ is the canonical basis of the algebra $\bA$ and $$\Aut_K(A)=\begin{cases}
\langle t_{\l_{m-4}}\rangle & \text{if char}(K)=0,\\
\langle t_{\l_{(m-4)_p}}\rangle& \text{if char}(K)=p>0,\\
\end{cases}
|\Aut_K(A)|=\begin{cases}
m-4 & \text{if char}(K)=0,\\
(m-4)_p& \text{if char}(K)=p>0,\\
\end{cases}$$ where $\l_{m-4}$ is a primitive $(m-4)$'th root of unity, $(m-4)_p$ is the $p$-co-prime divisor of $m-4$ if ${\rm char} (K)=p>0$.
\end{enumerate}
\end{enumerate}
\end{corollary}

{\it Proof}. 1. For polynomials of the type $x^a(1+\l_1x^{b_1}+\cdots +\l_ix^{b_i})$, where $1\leq b_1<\cdots <b_i$ and $\l_j\in K^\times$, that belong to the algebra $A$, $m\geq 2$ and $b_i\leq m-1-a\leq m-1-2=m-3$. If $a=2$ and $b_i\geq 1$ then $m\geq 4$ and $|\Aut_K(A)|\leq m-3$, by Theorem \ref{28Apr20}.(1). 

If $m$ is even and char$(K)=0$ then the algebra in the statement 2(a) has order $m-3$. So, the upper bound $m-3$ is sharp in this case.

If $m$ is even and char$(K)=p>0$ then 
$|\Aut_K(A)|\leq \max \{ i_p \, | \, 1\leq i \leq m-3 \}$, by Theorem \ref{28Apr20}.(1). 

If $m$ is odd then the number $m-3$ is even. So, for the algebra in the statement 2(a), the polynomials of the type $x^2(1+\l_1x^{b_1}+\cdots +\l_{m-3}x^{b_{m-3}})$ do not belong to its  canonical basis (since otherwise $m-1\in \G_A =\{ 2,4, \ldots , m-1\}$, a contradiction). Hence, $a\geq 3$. If $a=3$ then $1\leq b_i\leq m-1-a\leq m-1-3=m-4$, and so  $m\geq 5$ and $|\Aut_K(A)|\leq m-4$, by Theorem \ref{28Apr20}.(1).

If $m$ is odd and char$(K)=0$ then the algebra in the statement 2(b) has order $m-4$. So, the upper bound $m-4$ is sharp in this case.

If $m$ is odd and char$(K)=p>0$ then 
$\Aut_K(A)|\leq \max \{ i_p \, | \, 1\leq i \leq m-4 \}$, by Theorem \ref{28Apr20}.(1).

Suppose that  char$(K)=p$ and we fix $i$ such that $i_p=u(p,m)$, i.e. $1\leq i\leq m-3$, $i+j\leq m-1$ for some natural number $j=j(i)\geq 2$ such that $j\nmid m-1$. Let $g=x^j(1+x^i)$. Then, by the choice of the natural numbers $i$ and $j$,  $$A(g):=\sum_{k\geq 0}Kg^i+K[x]x^m\in \CA (m).$$
By Theorem \ref{28Apr20}.(3b), $|\Aut_K(A(g))|=i_p$. By Theorem \ref{28Apr20}.(3b), $u(p,m)\geq \max\{ |\Aut_K(A)|\; | \; A\in \CA (m), |\Aut_K(A)|<\infty\}$, and so the inequality is the equality.

2(a). Notice that $g^i=x^{2i}+\cdots $ for all $i\geq 1$, and so $\G_A= \{ 2i\, | \, 1\leq i\leq \frac{m-1}{2}\}$. Since $$g^i\equiv x^{2i} \mod (x^m)\;\; {\rm for \; all }\;\; i\geq 2,$$ the set $\{ 1, g, x^{2i}\, | \, 2\leq i \leq \frac{m-1}{2}\}$ is the canonical basis of the algebra $\bA$, and the rest  of the statement 2(a) follows from Theorem \ref{28Apr20}.(3). 

2(b). Similarly, $g^i=x^{3i}+\cdots $ for all $i\geq 1$, and so $\G_A= \{ 3i\, | \, 1\leq i\leq \frac{m-1}{3}\}$. Since $3\nmid m-1$, $m-1\not\in \G_A$. For all $i\geq 2$, $$g^i\equiv x^{3i} \mod (x^m).$$  So, the set  $\{ 1, g, x^{3i}\, | \, 2\leq i \leq \frac{m-1}{3}\}$ is the canonical basis of the algebra $\bA$, and the rest  of the statement 2(a) follows from Theorem \ref{28Apr20}.(3). 
 $\Box $\\
 
{\bf An explicit description of the set $\mO (m)$.}  For each element $\emptyset \neq \G \in \mS (m)$ and each prime number $p$, let
\begin{eqnarray*}
\mL (m,\G )&:=&\mL (\G ):=\{ l \, | \, 1\leq l \leq m-1, \;l+\G \not\subseteq \G \cup [m, \infty )\}= \{ l \, | \, 1\leq l \leq m-1, \; (l+\G )\cap C\G\neq \emptyset\}. \\
\mL_p (m,\G )&:=&\mL_p (\G ):=\{ l \, | \, 1\leq l \leq m-1,\;  p\nmid l, \; l+\G \not\subseteq \G \cup [m, \infty )\}=\{ l \, | \, 1\leq l \leq m-1,\;  p\nmid l, \\
& & \;(l+\G )\cap C\G\neq \emptyset\}.
\end{eqnarray*}
So, a natural number $l$ belongs to the set $\mL (m , \G)$ (resp., $\mL_p (m , \G)$) if 
  $1\leq l \leq m-1$ (resp., and $p\nmid l$) and $ l+\g \not\in\G \cup [m, \infty )$ for some element $\g \in \G$. If $l\in \mL (\G )$ then $$l\leq m-3$$ (since for $l\geq m-2$, $l+\G \subseteq [ m , \infty )$). By definition, $\mL (\emptyset ):=\emptyset$  and   $\mL_p (\emptyset ):=\emptyset$.
Let
$$
\mL (m,\G )_p:=\mL (\G )_p:=\{ l_p \, | \, l\in \mL (m, \G ) )\} $$
where $l_p$ is the $p$-co-prime divisor of $l$. In fact,
\begin{equation}\label{mpL=mLp}
\mL_p (m,\G )=\mL (m,\G )_p.
\end{equation}

{\it Proof.} By the very definition, $\mL_p (m,\G )\subseteq \mL (m,\G )$. Hence, $\mL_p (m,\G )\subseteq \mL (m,\G )_p$.

Suppose that  $l\in \mL (m,\G )$. We have to show that $l_p\in \mL_p (m,\G )$. Suppose that this is not true,  i.e. $l_p+\G \subseteq \G\cup [m, \infty)$. We seek a contradiction.  Recall that  $l=p^sl_p$ for some natural number $s\geq 0$. Then $l+\G = p^sl_p+\G \subseteq \G\cup [m, \infty)$, a contradiction. $\Box$

 Let 
\begin{eqnarray*}
 \mO (m) &:=& \{ |\Aut_K(A)|\; | \; A\in \CA (m),\;  |\Aut_K(A)|<\infty \},  \\
\mO (m, \G ) &:=& \{ |\Aut_K(A)|\; | \; A\in \CA (m, \G ),\;  |\Aut_K(A)|<\infty \}. 
\end{eqnarray*}
Clearly, $ \mO (m) =\bigcup_{\G \in \mS (m)}\mO (m, \G )$. Let 
$$ \mL (m) :=\bigcup_{\G \in \mS (m)}\mL (m , \G )\;\; {\rm and}\;\; \mL_p (m) :=\bigcup_{\G \in \mS (m)}\mL_p (m , \G )\stackrel{(\ref{mpL=mLp})}{=}\bigcup_{\G \in \mS (m)}\mL (m , \G )_p$$
where $p$ is a prime number. Theorem \ref{2May20} is an explicit description of the set $\mO (m)$.

\begin{theorem}\label{2May20}
Suppose that $m\geq 4$. Then 
$$\mO (m) = \begin{cases}
\mL (m)& \text{if char}(K)=0,\\
\mL_p(m)& \text{if char}(K)=p>0.\\
\end{cases}$$
\end{theorem}

{\it Remark.} Notice that the LHS  of the equality above depends on {\em all algebras} in $\CA (m)$ but the RHS depends only on {\em all semigroups} in $\mS (m)$ (which is a purely  combinatorial discrete object). \\

{\it Proof}. Let $R$ be the RHS of the equality in the theorem. 

(i) {\em For all $\G\in \mS (m)$, $\mO (m,\G ) \subseteq \begin{cases}
\mL (m, \G)& \text{if char}(K)=0,\\
\mL_p(m, \G)& \text{if char}(K)=p>0.\\
\end{cases}$ 

In particular,} $\mO (m) \subseteq R$:

 Given an algebra $A\in \CA (m)$ with finite automorphism group $G=\Aut_K(A)$.  By Theorem \ref{28Apr20}.(3), 
$$ |G|=\begin{cases}
\gcd (\bA ) & \text{if char}(K)=0,\\
\gcd_p (\bA )& \text{if char}(K)=p>0.\\
\end{cases}
$$
Let $\G = \G_A$ and $\{ 1, f_\g = x^\g u_\g\, | \, \g \in \G \}$ be the canonical basis of the algebra $\bA$ where $u_\g =1+\sum_{\d \in C\G (\g )} \l_{\g \d}x^{\d -\g}$ and $\l_{\g \d}\in K$. Since $l:=|G|<\infty$, $\l_{\g \d}\neq 0$ for some $\d\in C\G (\g )$ (by Corollary \ref{a2May20}). Hence, $\d -\g =\d'l$ for some natural number $\d'$ and 
$$ \g+ \d'l =\g +(\d -\g ) = \d \not\in \G \cup [m, \infty )$$
since $\d \in C\G (\g )$. Hence, 
 $$ l\in\begin{cases}
\mL (m, \G ) & \text{if char}(K)=0,\\
\mL_p (m, \G )& \text{if char}(K)=p>0,\\
\end{cases}
$$
 since otherwise $\G +l\subseteq \G \cup [m, \infty )$, and so $\G +il\subseteq \G \cup [m, \infty )$ for all natural numbers $i\geq 1$. In particular,   $\g +\d'l\in  \G \cup [m, \infty )$, a contradiction, and so the statement (i) follows.

(ii) $\mO (m) \supseteq R$: The inclusion follows from Proposition \ref{A2May20}.(2) (resp., Proposition \ref{A2May20}.(3)) if char$(K)=0$ (resp., char$(K)=p>0$). $\Box $\\

{\em Question. Is }
$$\mO (m, \G) = \begin{cases}
\mL (m, \G)& \text{if char}(K)=0,\\
\mL_p(m, \G)& \text{if char}(K)=p>0,\\
\end{cases}$$
{\em for each $\G \in \mS (m)$? (by the statement (i) in the proof of  Theorem \ref{2May20}, the LHS $\subseteq $ the RHS)}.

\begin{proposition}\label{A2May20}
Given $\emptyset \neq \G \in \mS (m)$ and $l\in \mL (\G )$. Then there exists an element $\g \in \G$ such that $\g +l\not\in \G$ and $\g +l\leq m-1$. Let $g=x^\g (1+x^l)$ and $A_{\g l}:= K+\sum_{i\geq 1}Kg^i+x^mK[x]$. Then 
\begin{enumerate}
\item $A_{\g l}\in \CA (m)$, $\G_{A_{\g l}}=\{i\g \, | \, 1\leq i\leq \frac{m-1}{\g}\}\subseteq \G$, $\ind (\G_{A_{\g l}})=\{\g \}$ and the polynomial $g$ is an element of the canonical basis of the algebra $A_{\g l}$.

\item $\Aut_K(A_{\g l})=\begin{cases}
\langle t_{\l_l}\rangle& \text{if char}(K)=0,\\
\langle t_{\l_{l_p}}\rangle& \text{if char}(K)=p>0.\\
\end{cases}$
\item Suppose that char$(K)=p>0$.  If, in addition, $p\nmid l$ then  $\Aut_K(A_{\g l})=\langle t_{\l_l}\rangle$.
\end{enumerate}
\end{proposition}

{\it Proof}. 1. Clearly, $A_{\g l}:= K\oplus\bigoplus_{i\geq 1}Kg^i\oplus x^mK[x]$. Hence, $A_{\g l}\in \CA (m)$, $\G_{A_{\g l}}=\{i\g \, | \, 1\leq i\leq \frac{m-1}{\g}\}\subseteq \G$ and  $\ind (\G_{A_{\g l}})=\{\g \}$. Since $\g +\l \in C\G (\g )$ and $\ind (\G_{A_{\g l}})=\{\g \}$,  the polynomial $g$ is an element of the canonical basis of the algebra $A_{\g l}$.


 
 

2. Statement 2 follows from statement 1 and Theorem \ref{28Apr20}.(3a,b).

3. Statement 3 is a particular case of statement 2. $\Box $\\
 
{ \em Definition.} Let $\mB (m):=\{ l \, | \, 1\leq l \leq m-1, l+\G\subseteq \G \cup [m,\infty )$ for all $\G \in \mS (m)\}$ where $m\geq 4$. \\
 
 By the very definition, the set $\mB (m)\cup \{ \infty\}$ is a subsemigroup of $\N_+/(m+\N )$, and 
\begin{equation}\label{LmBm}
\{ 1, \ldots , m-1\} = \mL (m) \coprod \mB (m)
\end{equation}
 is a disjoint union. Clearly, $$m-2, m-1\in \mB (m).$$ In view of Theorem \ref{2May20}, in order to find the set $\mL (m)$ it is much more faster, first, to find the set $\mB (m)$ and then apply (\ref{LmBm}). In general, $\mB (m) \neq \{ m-2, m-1\}$ (Lemma \ref{b3May20}.(1)) but the equality is also possible (Lemma \ref{b3May20}.(2)).\\
 
 {\em Question. Find a formula/lower bound/upper bound for the number $|\mB (m)|$.}

\begin{corollary}\label{c3May20}
Suppose that $m\geq 4$. Then 
$$\mO (m) = \begin{cases}
\{ 1, \ldots , m-1\}\backslash \mB (m)& \text{if char}(K)=0,\\
\{ l\, | \, l\in \{1, \ldots , m-1; p\nmid l\}\backslash \mB (m)\}=\{ l_p\, | \, l\in \{1, \ldots , m-1\}\backslash \mB (m)\} & \text{if char}(K)=p>0.\\
\end{cases}$$
\end{corollary}

{\it Proof}. The second equality is obvious.  Then   the corollary follows from Theorem \ref{2May20}, (\ref{mpL=mLp})  and (\ref{LmBm}). $\Box$

\begin{lemma}\label{a3May20}

\begin{enumerate}
\item For all $m\geq 4$, $1\not\in \mB (m)$. 
\item For each  $m\geq 4$, there is an algebra $A\in \CA (m)$ with $|\Aut_K(A)|=1$ (eg, $A=K\oplus Kx^{m-2}(1+x)\oplus x^mK[x]\in \CA (m)$).
\end{enumerate}
\end{lemma}

{\it Proof}. 1. Since $m\geq 4$, $\G :=\{ m-2\} \in \mS (m)$ and $1+\G = \{ m-1\} \not\subseteq \G\cup [m, \infty )$, i.e. $1\not\in \mB (m)$. 

2. Statement 2 follows from statement 1, Theorem \ref{2May20} and (\ref{LmBm}) (the fact that $|\Aut_K(A)|=1$ follows from Theorem \ref{28Apr20}.(3)). $\Box $

\begin{lemma}\label{b3May20}

\begin{enumerate}
\item Let $m=n!+1$ where $n\geq 3$. Then $\{ m-1-i\, | \, i=0,1,\ldots , n\}\subseteq \mB (m)$. In particular, $|\mB (m)|\geq n+1$.
\item Let $m=p+1$ where $p\geq 3$ is a prime number. Then $\mB (m) =\{ p-1, p\}$. 
\end{enumerate}
\end{lemma}

{\it Proof}. 1. The numbers $2,3,\ldots , n$ are divisors of $m-1=n!$. Therefore, for all $\emptyset \neq \G \in \mS (m)$, $\min (\G ) \geq n+1$. So, for all $i=0,1,\ldots , n$ and $\g \in \G$, 
$$m-1-i+\g \geq m-1-i+n+1=m+(n-i)\geq m,$$i.e. $m-1-i\in \mB (m)$. 

2. Since $m-1=p$ is a prime  number, $\G_i:=i\N_+\cap [2,p]\in \mS (m)$ for $i=2, \ldots , p-1$. Since $i\in \G_i$ and $i+(p-i)=p\not\in \G_i$, we have that $p-i\in \mL (m, \G_i)$. So, the elements $1, \ldots , p-2$ do not belong to the set $\mB (m)$. Since 
 $\{ p-1, p\} \subseteq \mB (m)$, we must have $\mB (m) =\{ p-1, p\}$. $\Box $\\
 
For a natural number $m\geq 2$ and a polynomial $p\in K[x]$, we denote by $p_{<m}$ a unique polynomial of degree $<m$ such that $p\equiv p_{<m}\mod (x^m)$.  

\begin{lemma}\label{a23Apr20}
Suppose that $m\geq 4$. Let  $g_l=x^{m-1-l}(1+x^l)=x^{m-1-l}+x^{m-1}$ where $1\leq l\leq m-3$ and $A_l=K+\sum_{i\geq 1}Kg_l^i+x^mK[x]$. Then $A_l \in \CA (m)$ and $A_l\neq K+x^mK[x]$ iff $ m-1-l\nmid m-1$. If the above equivalent conditions hold then 
\begin{enumerate}
\item $A_l=K\oplus Kg_l\oplus \bigoplus_{2\leq i<\frac{m-1}{m-1-l}}K(g_l^i)_{<m}\oplus x^mK[x]$ and $\G_{A_l}=\{ i(m-1-l)\, | \, 1\leq i<\frac{m-1}{m-1-l}\}$.
\item $\{ 1, g_l, x^{i(m-1-l)}\, | \, 2\leq i<\frac{m-1}{m-1-l}\}$ is the canonical basis of the algebra $A_l$.
\item  $\Aut_K(A_l)=\begin{cases}
\langle t_{\l_l}\rangle& \text{if char}(K)=0,\\
\langle t_{\l_{l_p}}\rangle& \text{if char}(K)=p>0.\\
\end{cases}$
\end{enumerate}
\end{lemma}

{\it Proof}. Since $g_l=x^{m-1-l}+x^{m-1}$ the `iff'-statement is obvious. 

1. Statement 1 follows from the fact that $g_l^i\equiv x^{i(m-1-l)}\mod (x^m)$ for all $i\geq 2$.

2. Statement 2 follows from statement 1.

3. Statement 3 follows from statement 2 and Theorem \ref{28Apr20}.(3). $\Box $\\

Clearly, $1\leq l\leq m-3$ $\Leftrightarrow $ $2\leq m-1-l<m-1$.
 The conditions
 $2\leq m-1-l<m-1$ and $ m-1-l| m-1$  (see Lemma \ref{a23Apr20}) are equivalent to the conditions: 
 $$ m-1-l=\frac{m-1}{i},\;\; 2\leq i \leq \frac{m-1}{2}\;\; {\rm and}\;\; i|m-1$$
 which are equivalent to the conditions:
$$ l=l(i), \;\;  2\leq i \leq \frac{m-1}{2}\;\; {\rm and}\;\; i|m-1$$
where $ l(i):=\frac{i-1}{i}(m-1)$. By Lemma \ref{a23Apr20}, 
\begin{equation}\label{Lm1}
\mL (m) \supseteq \begin{cases}
 \{ 1, \ldots , m-3\} \backslash \Big\{ l(i)\, | \,  2\leq i \leq \frac{m-1}{2}, \; i|m-1\Big\} &\text{if char}(K)=0,\\
 \{ l\, | \,  1\leq l\leq  m-3, p\nmid l\} \backslash \Big\{ l(i)_p\, | \,  2\leq i \leq \frac{m-1}{2}, \; i|m-1\Big\} &\text{if char}(K)=p>0.\\
\end{cases}
\end{equation}

Lemma \ref{d3May20} shows that under certain (explicit) conditions the elements $l(i)=\frac{i-1}{i}(m-1)$ in (\ref{Lm1}) can belong to the set $\mL (m)$.

\begin{lemma}\label{d3May20}
Suppose that $m\geq 4$, $i|m-1$ for some natural number $i$ such that $2\leq i \leq \frac{m-1}{2}$ and the number $n(i):=\frac{m-1}{i}-1\geq 2$ is not a divisor of the numbers $m-1$ and $m-2$. Let  $f_i=x^{n(i)}(1+x^{l(i)})$ where $l(i)=\frac{i-1}{i}(m-1)$, and $A_i=K+\sum_{j\geq 1}Kf_i^j+x^mK[x]$. Then 
\begin{enumerate}
\item $A_i\in \CA (m)$, $\G_{A_i}=\{ jn(i)\, | \, 1\leq j<\frac{m-1}{n(i)}\}$ and  $\{ 1, f_i, x^{jn(i)}\, | \, 2\leq j<\frac{m-1}{n(i)}\}$ is the canonical basis of the algebra $\bA_i$, and $f_i=x^{n(i)}+x^{m-2}$.
\item  $\Aut_K(A_i)=\begin{cases}
\langle t_{\l_{l(i)}}\rangle& \text{if char}(K)=0,\\
\langle t_{\l_{l(i)_p}}\rangle& \text{if char}(K)=p>0.\\
\end{cases}$
\end{enumerate}
\end{lemma}

{\it Proof}. Since $f_i^j=x^{jn(i)}+\cdots $ for all $j\geq 1$, $n(i)\geq 2$  and $n(i)\nmid m-1$, $$\G_{A_i}=\bigg\{ jn(i)\, | \, 1\leq j<\frac{m-1}{n(i)}\bigg\}\;\; {\rm and}\;\; A_i\in \CA (m).$$ Clearly,  $f_i=x^{n(i)}+x^{m-2}$. By the assumption, the number $m-2$ is not divisible by $n(i)$,  hence $m-2\not\in \G_{A_i}$. For all $j\geq 2$,
$$ f_i^j\equiv x^{jn(i)}\mod (x^m).$$
Therefore, the set $\{ 1, f_i, x^{jn(i)}\, | \, 2\leq j<\frac{m-1}{n(i)}\}$ is the canonical basis of the algebra $\bA_i$.

2. By statement 1, the set $\{ 1, f_i=x^{n(i)}(1+x^{l(i)}), x^{jn(i)}\, | \, 2\leq j<\frac{m-1}{n(i)}\}$ is the canonical basis of the algebra $\bA_i$. Hence, $\gcd (A)=l(i)$, and statement 2 follows from Theorem \ref{28Apr20}.(3). $\Box $\\

By (\ref{Lm1}) and Lemma \ref{d3May20}, if char$(K)=0$ then 
\begin{equation}\label{Lm2}
\mL (m) \supseteq \Bigg\{ l(i)\, | \,  2\leq i \leq \frac{1}{2}(m-1), \; i|m-1, \; n(i)\geq 2,\; \; n(i)\nmid m-1\; {\rm and}\; n(i)\nmid m-2\Bigg\}
\end{equation}
where $l(i)=\frac{i-1}{i}(m-1)$ and  $n(i)=\frac{m-1}{i}-1$; and if char$(K)=p>0$ then 
\begin{equation}\label{Lm3}
\mL (m) \supseteq \Bigg\{ l(i)_p\, | \,  2\leq i \leq \frac{1}{2}(m-1), \; i|m-1, \; n(i)\geq 2,\; n(i)\nmid m-1\; {\rm and}\; n(i)\nmid m-2\Bigg\}.
\end{equation}



\section{Generators and defining relations for the algebra  $\OO (\CA (m, \G ))$ of regular functions on the algebraic variety $\CA (m , \G )$}\label{MODULISP}

{\em In this section, the field $K$ is an arbitrary field, i.e. not necessarily algebraically closed (unless it is not stated otherwise)}.\\

The aim of the section is to  prove Lemma \ref{a21Jun20},  Theorem \ref{XX12May20} and  Theorem \ref{XY12May20} which show that the set $\CA (m, \G )$ is an affine algebraic variety and  gives  explicit sets of  generators and defining relations for the algebra $\OO (\CA (m, \G ))$ of regular functions on $\CA (m, \G )$.\\

{\bf The action of the algebraic torus $\mT$ on $\CA (m)$.}  The algebraic torus $\mT =\{ t_\l \, | \, \l \in K^\times \}$ acts on the set $\CA (m)$ by the rule:
$$\mT \times \CA (m)\ra \CA (m), \;\; (t_\l , A)\mapsto t_\l (A).$$
The subsets $\CA (m, \G )$ of $\CA (m)$, where $\G \in \mS (m)$, are $\mT$-stable (that is $\mT \CA (m, \G ) =\CA (m , \G )$). By Theorem \ref{3May20}.(1), {\em the set of $\mT$-orbits of} $\CA (m)$,
\begin{equation}\label{AAm0}
\mA (m):= \CA (m)/\mT = \{ \mT A\, | \, A\in \CA (m)\},
\end{equation}
{\em is the set of isomorphism classes of algebras } $\CA (m)$. By (\ref{grbA2}), 
\begin{equation}\label{AAm}
\mA (m)=\coprod_{\G \in \mS (m)}\mA (m , \G )\;\; {\rm where}\;\; \mA (m,\G ) := \CA (m, \G )/\mT
\end{equation}
is the set of isomorphism classes of algebras in $\CA (m, \G )$. \\

{\bf The set $\CA (m, \G )$ is an affine algebraic variety.}  Let $A\in \CA (m , \G )$. 
  Recall that  $\ind (\G )=\{ \nu_1, \ldots , \nu_s\}$ and $f_{\nu_i}=x^{\nu_i}+\sum_{j\in C\G (\nu_i)}\l_{\nu_i, j}x^j$ for $i=1, \ldots , s$.
The coefficients $\{ \l_{\nu_i, j}\, | \, i=1, \ldots , s; j\in C\G (\nu_i)\}$ uniquely determined the algebra $A$. We treat them as regular functions on the algebraic variety $\CA (m, \G )$.

\begin{lemma}\label{a21Jun20}

\begin{enumerate}
\item  Suppose that $|\ind (\G )|=1$, i.e. $\ind (\G )=\{ \nu_1\}$. Then the affine algebraic variety $\CA (m, \G )$ is isomorphic to the affine space $\mA^{|C\G (\nu_1)|}$.
\item  Suppose that $s=|\ind (\G )|\geq 2$ and $\dec (\G)_{\geq 2}=\emptyset$. Then the affine algebraic variety $\CA (m, \G )$ is isomorphic to the affine space $\mA^{|C\G (\nu_1)|+\cdots +|C\G (\nu_s)|}$ where $\ind (\G )=\{ \nu_1, \ldots , \nu_s\}$.
\end{enumerate}
\end{lemma}

{\it Proof}. 1. Statement 1 follows from Lemma \ref{Aa7May20}.

2. Statement 2 follows from Corollary \ref{Ab7May20}. $\Box $\\

So, it remains to consider the case when $s=|\ind (\G )|\geq 2$ and $\dec (\G)_{\geq 2}\neq \emptyset$. The algebra of differential operators $\CD (K[x])$ on the polynomial algebra $K[x]$  is equal to 
$$ \CD (K[x])=\bigoplus_{n\geq 0}K[x]\der^{[n]}\;\; {\rm where}\;\; \der^{[n]}:=\frac{\der^n}{n!}, \;\;\der := \frac{d}{dx}$$ and the action of the differential operators $\der^{[n]}$ on the polynomial algebra $K[x]$ is given by the rule
$$\der^{[n]}(x^m)=\begin{cases}
{m\choose n}x^{m-n}& \text{if }m\geq n,\\
0& \text{if }m<n.\\
\end{cases}
$$
If the field $K$ has characteristic zero then the algebra $\CD (K[x])$ is generated by the elements $x$ and $\der$ that satisfy the defining relation $\der x -x\der =1$, and the algebra $\CD (K[x])$ is called the {\em Weyl algebra}. 

Recall that for  a polynomial $f\in K[x]$,  we denoted by  $c_j(f)$ the coefficient of $x^j$. Clearly,
\begin{equation}\label{cjf}
c_j(f)=\der^{[j]}(f)|_{x=0}.
\end{equation}

Recall that $A\in \CA (m , \G )$;   $s=|\ind (\G )|\geq 2$ and $\dec (\G)_{\geq 2}\neq \emptyset$;   for each element $\g \in \G$, we fixed an element $a(\g )$, see (\ref{nuQ2}); $\ind (\G )=\{ \nu_1, \ldots , \nu_s\}$ and $f_{\nu_i}=x^{\nu_i}+\sum_{j\in C\G (\nu_i)}\l_{\nu_i, j}x^j$ for $i=1, \ldots , s$. The elements $\{ 1, f^{a(\g )}\, | \, \g \in \G\}$ is a $K$-basis of the algebra $\bA$. 
 For all $i=1, \ldots , s$ and $\g \in \G$, we have the equality in the algebra $\bA$, 
\begin{equation}\label{cjf1}
f_{\nu_i}f^{a(\g )}=f^{a(\nu_i +\g)}+\sum_{\d\in \G (\nu_i+\g )} \eta_{\nu_i, a(\g ); \d}f^{a(\d )},
\end{equation}
where for each $i=1, \ldots , s$, the coefficients 
$\eta_{\nu_i, a(\g ), \d}$, $\d \in \G$,  are the {\em unique solution} of the uni-triangular system of $|\G |$ linear equations 
(the $|\G |\times |\G |$ matrix of which has diagonal elements that are equal to 1):
\begin{equation}\label{cjf2}
 \eta_{\nu_i, a(\g ); \d}+\sum_{\d'\in \G (\nu_i+\g , \d)} \eta_{\nu_i, a(\g ); \d'}c_\d\Big(f^{a(\d' )}\Big)=c_\d \Big(f_{\nu_i}f^{a(\g )}-f^{a(\nu_i +\g)}\Big), \;\; \g \in \G . 
\end{equation}
Each  element $\eta_{\nu_i, a(\g ); \d}$, $\d \in \G$, is an {\em  explicit polynomial} in the variables $\{ \l_{\nu_j, \g'}\, | \, 1\leq j \leq s, \; \g'\in C\G (\nu_j), \; \g'\leq \d \}$.
 The elements $\eta_{\nu_i, a(\g ); \d}$, $\d \in \G$, can be found recursively using (\ref{cjf2}). 


If $a(\nu_i +\g ) = e_i+a(\g )$ then the equality (\ref{cjf1}) is the tautology $f_{\nu_i}f^{a(\g )}=f_{\nu_i}f^{a(\g )}$, i.e. all the coefficients  
 $\eta_{\nu_i, a(\g ); \d}$ are equal to zero, and vice versa. In particular, if $\nu_i+\g \not\in \dec_{\geq 2}(\G )$ then $a(\nu_i +\g ) = e_i+a(\g )$.

Theorem \ref{XX12May20} shows that the 
 set $\CA (m, \G )$ is an affine algebraic variety and gives an explicit  set of its defining equations (see also   Theorem \ref{XY12May20} for another approach). 
 
  \begin{theorem}\label{XX12May20}
We keep the notation as above (recall that for each element $\g \in \G$, we fix an element $a(\g )$, see (\ref{nuQ2})).  Suppose that $s=|\ind (\G )|\geq 2$ and $\dec (\G)_{\geq 2}\neq \emptyset$.
 Then the  set $\CA (m, \G )$ is an affine algebraic variety and  the algebra of regular functions on it, $\OO (\CA (m, \G ))$, is a factor algebra of the polynomial algebra $P(m,\G ):=K[\l_{\nu_i, j}\,| \,  i=1, \ldots , s; j\in C\G (\nu_i)]$ in $n(m,\G ):=\sum_{i=1}^s |C\G (\nu_i)|$ variables $\l_{\nu_i, j}$ by the defining relations: For each pair $(\nu_i, \g )\in \ind (\G )\times \G$  such that $\nu_i +a(\g ) \neq a(\nu_i+\g )$  and  $j\in C\G (\nu_i+\g )$, 
 $$ c_j\Big(f_{\nu_i}f^{a(\g )}-f^{a(\nu_i +\g)}-\sum_{\d\in \G (\nu_i+\g )} \eta_{\nu_i, a(\g ); \d}f^{a(\d )}\Big)=0, \;\; {\rm i.e.}$$
$$ \der^{[j]}\Big(f_{\nu_i}f^{a(\g )}-f^{a(\nu_i +\g)}-\sum_{\d\in \G (\nu_i+\g )} \eta_{\nu_i, a(\g ); \d}f^{a(\d )}\Big)|_{x=0}=0$$  
where the polynomials $\eta_{\nu_i, a(\g ); \d}$ are defined in (\ref{cjf2}).  
\end{theorem}
 
 {\it Proof}. By the very definition,  the generators $\{ \l_{\nu_i, j}\}$ of the algebra $\OO (\CA (m, \G ))$ satisfy the relations of the theorem.
 
Conversely, suppose that the  {\em scalars} 
$\{ \l_{\nu_i, j}\}$ are a solution to the system of  equations of the theorem. They determine the elements $f_{\nu_i}=x^{\nu_i}+\sum_{j\in C\G (\nu_i)}\l_{\nu_i, j}x^j$, $i=1, \ldots , s$ of the algebra $K[x]/(x^m)$.  
We have to show that the subalgebra 
$\bA'$ of $K[x]/(x^m)$, which is generated by the elements $f_{\nu_1},  \ldots ,f_{\nu_s}$, is equal to $$V:=K\oplus\bigoplus_{\g \in \G} Kf^{a(\g )}$$ (since then the subalgebra $A'$ of $K[x]$, which is generated by the elements 
$f_{\nu_1},  \ldots ,f_{\nu_s}$ and the ideal $(x^m)$ of $K[x]$, would belong to $\CA (m, \G))$.

 Clearly, $V\subseteq \bA'$ (by the definition of the elements $f_{\nu_i}$). The defining relations of the theorem mean that the equality (\ref{cjf1}) holds for all pairs $(\nu_i, \g )\in \ind (\G )\times \G$ such that  $\nu_i +a(\g ) \neq a(\nu_i+\g )$.  But for all  pairs $(\nu_i, \g )\in \ind (\G ) \times \G  $ such that $a(\nu_i+\g ) = \nu_i+a(\g )$ the equality (\ref{cjf1})   holds automatically, it is  simply the tautology  $f_{\nu_i}f^{a(\g )}=f_{\nu_i}f^{a(\g )}$. So, the equality (\ref{cjf1})  holds for all elements $(\nu_i, \g )\in \ind (\G ) \times \G $, i.e. $\bA'=V$, as required.  $\Box$ \\

 Let $\bA'$ be a subalgebra of $K[x]/(x^m)$ which is generated by the elements $f_{\nu_i}=x^{\nu_i}+\sum_{j\in C\G (\nu_i)}\l_{\nu_i, j}x^j$, $i=1, \ldots , s$, where we treat the coefficients $\l_{\nu_i, j}$ as independent parameters (indeterminates).  For each $\g \in \dec_{\geq 2}(\G)$ and $b\in \Rel (\g )\backslash \{ a(\g )\}$, consider the equation (see Theorem \ref{A6May20}), 
$$ f^b = f^{a(\g)}+\sum_{\g'\in \G (\g )}\th_{\g , \g';b}f^{a(\g')}.
 $$
 The coefficients $\th_{\g , \g';b}$ are a {\em unique solution} to the uni-triangular system of  $|\G(\g )|$ linear equations (the $|\G(\g )|\times |\G(\g )|$ matrix of which has diagonal elements that are equal to 1): 
\begin{equation}\label{cjf3}
 \th_{\g , \g';b} +\sum_{\g''\in \G(\g, \g')}\th_{\g , \g'';b}\, c_{\g'}\Big(f^{a(\g'')}\Big)=c_{\g'}\Big(f^b-f^{a(\g)}\Big)
\end{equation}
 where $ \G(\g, \g'):=\{ \g''\in \G \, | \, \g <\g''<\g'\}$.  Each  element $\th_{\g , \g';b}$, $\g' \in \G (\g )$, is an {\em  explicit polynomial} in the variables $\{ \l_{\nu_j, \d}\, | \, 1\leq j \leq s, \; \d\in C\G (\nu_j), \; \d\leq \g'\}$. Theorem \ref{XY12May20} also shows that the 
 set $\CA (m, \G )$ is an affine algebraic variety and gives an explicit  set of its defining equations. It also gives the lower bound for the dimension of the algebraic variety $\CA (m, \G )$. 
 
\begin{theorem}\label{XY12May20}
We keep the notation as above.   Suppose that $s=|\ind (\G )|\geq 2$ and $\dec (\G)_{\geq 2}\neq \emptyset$.
  The set $\CA (m, \G )$ is an affine algebraic variety, the algebra of regular functions on it, $\OO (\CA (m, \G ))$, is a factor algebra of the polynomial algebra $P(m,\G )=K[\l_{\nu_i, j}\,| \,  i=1, \ldots , s; j\in C\G (\nu_i)]$ in $n(m,\G )=\sum_{i=1}^s |C\G (\nu_i)|$ variables $\l_{\nu_i, j}$ by the defining relations: For each $\g \in \dec_{\geq 2}(\G)$ and $b\in \Rel (\g )\backslash \{ a(\g )\}$,
$$ c_j(f^b )= c_j(f^{a(\g)})+\sum_{\g'\in \G (\g )}\th_{\g , \g';b}\, c_j(f^{a(\g')}), \;\; j\in C\G (\g ),
 $$
 where the polynomials  $\th_{\g , \g';b}$ are defined in (\ref{cjf3}).
 The number of equations is $$l(m, \G ):=\sum_{\g \in \dec_{\geq 2}(\G )} \Big(|\Rel (\g )|-1\Big)\cdot |C\G (\g )|.$$ In particular, the dimension of the variety $\CA (m, \G)$ is not smaller than $n(m, \G ) - l(m, \G )$. 
\end{theorem}

{\it Proof}.  By Theorem \ref{A6May20}, the generators $\{ \l_{\nu_i, j}\}$ of the algebra $\OO (\CA (m, \G ))$ satisfy the relations of the theorem.
 
Conversely, suppose that the  {\em scalars} 
$\{ \l_{\nu_i, j}\}$ are a solution to the equations of the theorem. They determine the elements $f_{\nu_i}=x^{\nu_i}+\sum_{j\in C\G (\nu_i)}\l_{\nu_i, j}x^j$, $i=1, \ldots , s$ of the algebra $K[x]/(x^m)$.  
We have to show that the subalgebra 
$\bA'$ of $K[x]/(x^m)$, which is generated by the elements $f_{\nu_1},  \ldots ,f_{\nu_s}$, is equal to $$V:=K\oplus\bigoplus_{\g \in \G} Kf^{a(\g )}$$ (since then the subalgebra $A'$ of $K[x]$, which is generated by the elements 
$f_{\nu_1},  \ldots ,f_{\nu_s}$ and the ideal $(x^m)$ of $K[x]$, would belong to $\CA (m, \G))$.

 Clearly, $V\subseteq \bA'$. On the other hand,  the defining relations of the theorem and (\ref{cjf3}) mean that for each $\g \in \dec_{\geq 2}(\G)$ and $b\in \Rel (\g )\backslash \{ a(\g )\}$,
 $$ f^b = f^{a(\g)}+\sum_{\g'\in \G (\g )}\th_{\g , \g';b}\, f^{a(\g')}.
 $$
 Now, by (\ref{BGbas1}),    $\bA' =K\oplus\bigoplus_{\g \in \G} Kf^{a(\g )}$. $\Box$

   \begin{corollary}\label{a27Jun20}
 The set of irreducible components of the affine algebraic variety $\CA (m)$ is the union of irreducible components of the affine algebraic varieties $\{ \CA (m , \G )\, | \, \G \in \mS(m)\}$. 
\end{corollary}

{\it Proof}.  The affine algebraic variety  $\CA (m) = \coprod_{\G \in \mS (m)}\CA (m, \G )$ is a finite disjoint union of its closed subsets $\CA (m, \G )$ (Theorem \ref{XX12May20} or Theorem \ref{XY12May20}) and the statement follows.   $\Box $\\

Recall that $f_{\nu_i}=x^{\nu_i}+\sum_{j\in C\G (\nu_i)}\l_{\nu_i, j}x^j$ for $i=1, \ldots , s$, and the  algebraic  group $\mT$ acts on the algebraic variety $\CA (m, \G )$. The  action of the algebraic group on the algebra 
 $\OO (\CA (m, \G ))$ is given by the rule: For all $i=1, \ldots , s$, $j\in C\G (\nu_i)$ and $\l \in K^\times$, 
\begin{equation}\label{AAm1}
t_\l (\l_{\nu_i, j})=\l^{-\nu_i+j}\l_{\nu_i, j}
\end{equation} 
since $t_\l (f_{\nu_i})=\l^{\nu_i}\Big(x^{\nu_i}+\sum_{j\in C\G (\nu_i)}\l^{-\nu_i+j}\l_{\nu_i, j}x^{j-\nu_i}\Big)$. \\

Recall that for each primitive $n$'th root of unity $\l_n$,  the group $\mT$ contains the cyclic group of order $n$, $C_n=\langle t_{\l_n}\rangle$. The groups $\{ C_n\}$ are all the finite subgroups of $\mT$. 

\begin{corollary}\label{a12May20}
($K$ is an algebraically closed field). Given  $n\in \mO (m, \G )$.
\begin{enumerate}
\item The set $\CA (m, \G )^{C_n}$ of fixed points of the group $C_n=\langle t_{\l_n}\rangle\subseteq \mT$ (where $\l_n$ is a primitive $n$'th root of unity) is a non-empty closed subvariety of the affine algebraic variety $\CA (m , \G )$ the defining equations of which are 
$\l_{\nu_i, j}=0$ for all $i=1, \ldots , s$ and $j\in C\G (\nu_i)$ such that $n\nmid j-\nu_i$. The set $\CA (m, \G )^{C_n}$ contains precisely all the algebras $A$ in $\CA (m, \G )$ with $C_n\subseteq \Aut_K(A)$. 
\item The set $\CA (m, \G )^{C_n}/\mT$ is a set of isomorphism classes of algebras $A$ in $\mA (m, \G )$ with $C_n\subseteq \Aut_K(A)$. If $n|n'$ then $\CA (m, \G )^{C_n}/\mT\supseteq \CA (m, \G )^{C_{n'}}/\mT$. 
\item The set $\CA (m, \G )^{C_n}/\mT\backslash \bigcup_{n\neq l\in \mO (m, \G ), n|l}\CA (m, \G )^{C_l}/\mT$ is a set of isomorphism classes of algebras $A$ in $\mA (m, \G )$ with $C_n= \Aut_K(A)$.  
\item The set $\CA (m, \G )/\mT\backslash \Big\{ \{ A_{mon} (\G )\} \cup \bigcup_{1\neq l\in \mO (m, \G )}\CA (m, \G )^{C_l}/\mT\Big\}$ is a set of isomorphism classes of algebras $A$ in $\mA (m, \G )$ with $ \Aut_K(A)=\{ e\}$. It is an open set of the algebraic variety $\CA (m , \G )$.
\item $\CA (m, \G )^\mT=\{ A_{mon}(\G )\}$.
\end{enumerate}
\end{corollary}

{\it Proof}. 1. Statement 1 follows at once from 
 (\ref{AAm1}).
 
2.  Statement 2 follows from statement 1. 

3. Statement 3 follows from statement 2.

5. Statement 5 is obvious.

4. Statement 4 follows from statements 1, 3 and  5. $\Box $ \\

$${\bf Acnowledgements} $$

The author would like to thank the Royal Society for support.

Department of Pure Mathematics

University of Sheffield

Hicks Building

Sheffield S3 7RH

UK

email: v.bavula@sheffield.ac.uk

\end{document}